\documentclass[11pt]{article} 
\usepackage{latexsym} 
\usepackage{amsmath}
\usepackage{amssymb} 
\usepackage{amscd} 
\usepackage{graphicx}
\usepackage{layout}
\usepackage{pifont}

\begin{document} 

\newtheorem{Th}{Theorem}[section]
\newtheorem{Cor}{Corollary}[section]
\newtheorem{Prop}{Proposition}[section]
\newtheorem{Lem}{Lemma}[section]
\newtheorem{Def}{Definition}[section]
\newtheorem{Rem}{Remark}[section]
\newtheorem{Ex}{Example}[section]
\newtheorem{stw}{Proposition}[section]


\newcommand{\bet}{\begin{Th}}
\newcommand{\ent}{\stepcounter{Cor}
   \stepcounter{Prop}\stepcounter{Lem}\stepcounter{Def}
   \stepcounter{Rem}\stepcounter{Ex}\end{Th}}


\newcommand{\bec}{\begin{Cor}}
\newcommand{\enc}{\stepcounter{Th}
   \stepcounter{Prop}\stepcounter{Lem}\stepcounter{Def}
   \stepcounter{Rem}\stepcounter{Ex}\end{Cor}}
\newcommand{\bep}{\begin{Prop}}
\newcommand{\enp}{\stepcounter{Th}
   \stepcounter{Cor}\stepcounter{Lem}\stepcounter{Def}
   \stepcounter{Rem}\stepcounter{Ex}\end{Prop}}
\newcommand{\bel}{\begin{Lem}}
\newcommand{\enl}{\stepcounter{Th}
   \stepcounter{Cor}\stepcounter{Prop}\stepcounter{Def}
   \stepcounter{Rem}\stepcounter{Ex}\end{Lem}}
\newcommand{\bef}{\begin{Def}}
\newcommand{\enf}{\stepcounter{Th}
   \stepcounter{Cor}\stepcounter{Prop}\stepcounter{Lem}
   \stepcounter{Rem}\stepcounter{Ex}\end{Def}}
\newcommand{\ber}{\begin{Rem}}
\newcommand{\enr}{
   \stepcounter{Th}\stepcounter{Cor}\stepcounter{Prop}
   \stepcounter{Lem}\stepcounter{Def}\stepcounter{Ex}\end{Rem}}
\newcommand{\bee}{\begin{Ex}}
\newcommand{\ene}{
   \stepcounter{Th}\stepcounter{Cor}\stepcounter{Prop}
   \stepcounter{Lem}\stepcounter{Def}\stepcounter{Rem}\end{Ex}}
\newcommand{\Proof}{\noindent{\it Proof\,}:\ }

\newcommand{\EE}{\mathbf{E}}
\newcommand{\QQ}{\mathbf{Q}}
\newcommand{\R}{\mathbf{R}}
\newcommand{\C}{\mathbf{C}}
\newcommand{\ZZ}{\mathbf{Z}}
\newcommand{\KK}{\mathbf{K}}
\newcommand{\NN}{\mathbf{N}}
\newcommand{\PP}{\mathbf{P}}
\newcommand{\HH}{\mathbf{H}}
\newcommand{\uuu}{\boldsymbol{u}}
\newcommand{\xxx}{\boldsymbol{x}}
\newcommand{\aaa}{\boldsymbol{a}}
\newcommand{\bbb}{\boldsymbol{b}}
\newcommand{\AAA}{\mathbf{A}}
\newcommand{\BBB}{\mathbf{B}}
\newcommand{\ccc}{\boldsymbol{c}}
\newcommand{\iii}{\boldsymbol{i}}
\newcommand{\jjj}{\boldsymbol{j}}
\newcommand{\kkk}{\boldsymbol{k}}
\newcommand{\rrr}{\boldsymbol{r}}
\newcommand{\FFF}{\boldsymbol{F}}
\newcommand{\yyy}{\boldsymbol{y}}
\newcommand{\ppp}{\boldsymbol{p}}
\newcommand{\qqq}{\boldsymbol{q}}
\newcommand{\nnn}{\boldsymbol{n}}
\newcommand{\vvv}{\boldsymbol{v}}
\newcommand{\eee}{\boldsymbol{e}}
\newcommand{\fff}{\boldsymbol{f}}
\newcommand{\www}{\boldsymbol{w}}
\newcommand{\0}{\boldsymbol{0}}
\newcommand{\lon}{\longrightarrow}
\newcommand{\ga}{\gamma}
\newcommand{\pa}{\partial}
\newcommand{\QED}{\hfill $\Box$}
\newcommand{\id}{{\mbox {\rm id}}}
\newcommand{\Ker}{{\mbox {\rm Ker}}}
\newcommand{\grad}{{\mbox {\rm grad}}}
\newcommand{\ind}{{\mbox {\rm ind}}}
\newcommand{\rot}{{\mbox {\rm rot}}}
\newcommand{\diver}{{\mbox {\rm div}}}
\newcommand{\Gr}{{\mbox {\rm Gr}}}
\newcommand{\LG}{{\mbox {\rm LG}}}
\newcommand{\Diff}{{\mbox {\rm Diff}}}
\newcommand{\Symp}{{\mbox {\rm Symp}}}
\newcommand{\Ct}{{\mbox {\rm Ct}}}
\newcommand{\Uns}{{\mbox {\rm Uns}}}
\newcommand{\rank}{{\mbox {\rm rank}}}
\newcommand{\sign}{{\mbox {\rm sign}}}
\newcommand{\Spin}{{\mbox {\rm Spin}}}
\newcommand{\Sp}{{\mbox {\rm sp}}}
\newcommand{\Int}{{\mbox {\rm Int}}}
\newcommand{\Hom}{{\mbox {\rm Hom}}}
\newcommand{\codim}{{\mbox {\rm codim}}}
\newcommand{\ord}{{\mbox {\rm ord}}}
\newcommand{\Iso}{{\mbox {\rm Iso}}}
\newcommand{\corank}{{\mbox {\rm corank}}}
\def\mod{{\mbox {\rm mod}}}
\newcommand{\pt}{{\mbox {\rm pt}}}
\newcommand{\enP}{\hfill $\Box$ \par\vspace{5truemm}}
\newcommand{\qed}{\hfill $\Box$ \par}
\newcommand{\spe}{\vspace{0.4truecm}}
\newcommand{\no}{\noindent}
\newcommand{\Tan}{{\mbox{\rm Tan}}}
\newcommand{\Reg}{{\mbox {\rm Reg}}}
\newcommand{\Lag}{{\mbox {\rm {\scriptsize Lag}}}}
\newcommand{\GL}{{\mbox {\rm GL}}}
\newcommand{\IGr}{{\mbox {\rm IGr}}}
\newcommand{\type}{{\mbox {\rm type}}}
\newcommand{\Jetcod}{{\mbox {\rm Jet-codim}}}

\newcommand{\dint}[2]{{\displaystyle\int}_{{\hspace{-1.9truemm}}{#1}}^{#2}}


\title{
Singularities of Tangent Varieties \\
to Curves and Surfaces 
}

\author{Goo Ishikawa
\\
Department of Mathematics, Hokkaido University, Japan
\\
e-mail: ishikawa@math.sci.hokudai.ac.jp
}

\date{ }

\maketitle

\begin{center}
Abstract. 
\end{center}
It is given the diffeomorphism classification on generic singularities of tangent varieties to curves 
with arbitrary codimension in a projective space. 
The generic classifications are performed in terms of certain geometric structures and 
differential systems on flag manifolds, via several techniques in differentiable algebra. 
It is provided also the generic diffeomorphism classification of singularities on tangent varieties 
to contact-integral curves in the standard contact projective space. 
Moreover we give basic results on the classification of  
singularities of tangent varieties to generic surfaces and Legendre surfaces. 

\section{Introduction}

Embedded tangent spaces to a submanifold draw a variety in the ambient space, 
which is called the {\it tangent variety} to the submanifold. 
Tangent varieties appear in various geometric problems and applications naturally. 
See for instance \cite{AG}\cite{FP}\cite{BG}. 
Developable surfaces, varieties with degenerate Gauss mapping and 
varieties with degenerate projective dual are obtained by tangent varieties. 
Tangent varieties provide several important examples of non-isolated singularities in applications of geometry. 
We observe relations of tangent varieties to invariant theory and geometric theory of differential equations 
(see \cite{LW}, also see 
Examples \ref{umbilical bracelet} and \ref{Tangents to Veronese}). 

It is known, in the three dimensional Euclidean space, 
that the tangent variety (tangent developable) to a generic space curve 
has singularities each of which is locally diffeomorphic to the {\it cuspidal edge} 
or to the {\it folded umbrella} ({\it cuspidal cross cap}), as is found by Cayley and Cleave \cite{Cleave}. 
Cuspidal edge singularities appear along ordinary points, while 
the folded umbrella appears at an isolated point of zero torsion \cite{BG}\cite{Porteous}. 

The classification was generalised to more degenerate cases 
by Mond \cite{Mond1}\cite{Mond2} and Scherbak \cite{Scherbak}\cite{Arnol'd3} 
and applied to various geometry (see for instance \cite{CI}\cite{IKY}). 
If we consider a curve together with its osculating framings, we are led to the classification of 
tangent varieties to generic osculating framed curves, possibly with singularities in themselves, 
in the three dimensional space. 
Then the list consists of 4 singularities: 
{\it cuspidal edge}, {\it folded umbrella} and 
moreover {\it swallowtail} and {\it Mond surface} (\lq cuspidal beak to beak\rq) \cite{Ishikawa12}. 
However the author could not find any literature treating the classification 
of singularities appearing in tangent varieties to higher codimensional curves. 

The diffeomorphism types of tangent varieties to curves are 
invariant under projective transformations. 
In this paper, we consider curves in projective spaces and 
show the classification results on generic singularities of tangent 
varieties to curves with arbitrary codimension in projective spaces. 

The tangent variety can be defined for a \lq frontal' variety. 
A frontal variety has the well-defined embedded tangent space at each point, 
even where the variety is singular. 
In Cauchy problem of single unknown function, 
we have wave-front sets, which are singular hypersurfaces \cite{Arnol'd2}. 
They are called fronts and form an important class of frontal varieties. 
Also higher codimensional wave-fronts are examples of frontal varieties, 
which appear in, for instance, Cauchy problem of several unknown functions, 
where initial submanifolds of arbitrary codimension evolve to frontal varieties (cf. \cite{Givental'}\cite{Kossowski1}). 

First, in \S \ref{Frontal maps and tangent varieties.},  
we introduce the notion of frontal maps and 
frontal varieties,  generalising that of submanifolds and fronts (Definition \ref{def-frontal}). 
Moreover we define their tangent maps and tangent varieties (Definition \ref{def-tangent}). 
Then we give the classification of tangent varieties to generic curves in projective spaces 
(Theorem \ref{tangent to generic curve}). 
In fact we find that 
the tangent variety to a generic curve in $\R P^{N+1}$ has the unique singularity,  
the higher codimensional cuspidal edge, if $N+1 \geq 4$. 

In the geometric theory of curves, however, we usually treat not just curves but 
we attach an appropriate frame with curves. 
Thus, to solve the generic classification problem properly, 
we relate the study of tangent varieties 
to certain kinds of differential systems on appropriate flag manifolds 
in \S \ref{Differential systems on flag manifolds.}. 
Note that the method was initiated by Arnol'd and Scherbak \cite{Scherbak}. 
Also note that it is standard to use flag manifolds in the theory of space curves (\cite{Wall2}). 
We can utilise various types of flag manifolds. 
In fact, in this paper,  we select three kinds of flag manifolds and three kinds of differential systems, 
correspondingly to the classes of curves endowed with osculating-frames, with tangent-frames and 
with tangent-principal-normal-frames. 
Then we present the classification results on the singularities of 
which generically appear for these three kinds of classes 
of curves in projective spaces 
(Theorems \ref{tangent-framed classification}, 
\ref{tangent-principal-normal-framed classification}, \ref{osculating-frame-geq4}).

In \S \ref{Flat projective structure and the type of a curve.}, the notion of types of curve-germs 
are recalled. Curves of finite type are frontal and their tangent varieties are frontal. 
We classify the generic types of curves, and then we show a kind of determinacy of 
the tangent variety for 
each generic type of curves. 

In \S \ref{Codimension formulae and transversality theorems.}, we classify the list of types of generic 
curves satisfying geometric conditions. 
To do this, we establish the codimension formulae giving the codimension of the set of curves, for given type, 
which satisfy a given geometric integrality condition in each case. 
Then the transversality theorem 
implies the restriction on types of generic curves.  

In \S \ref{Opening of differentiable map-germs.}, we introduce the key notion of openings of differentiable 
map-germs, which has close relations with that of frontal varieties. We collect necessary results on 
differentiable algebras to solve the generic classification problems treated in this paper. 
Moreover,  in \S \ref{Normal forms of tangent surfaces.}, 
using the method of differentiable algebra, 
we show the normal forms of tangent varieties appearing in the generic classification problems we have treated 
in this paper. In particular the main results in this paper, Theorems \ref{tangent to generic curve}, 
\ref{tangent-framed classification}, \ref{tangent-principal-normal-framed classification} 
and \ref{osculating-frame-geq4} are proved. 

In \S \ref{Singularities on tangent varieties to osculating framed contact-integral curves.}, 
we treat contact-integral curves and their tangent varieties. 
If $V$ is a symplectic vector space $V$, then the projective space 
$P(V)$ has the canonical contact structure. 
Then we give the generic diffeomorphism classification 
of singularities on tangent varieties to 
\lq osculating framed contact-integral' curves in $P(V)$ (Theorems  
\ref{osculating framed contact-integral-geq7}, \ref{osculating framed contact-integral-=5}). 
For this, in particular, 
we show that the diffeomorphism type of $\Tan(\gamma)$ is unique for a 
curve of type $(1, 3, 4, 6)$ in $\R P^4$ in this paper. 
Note that it is known that the diffeomorphism type of $\Tan\Tan(\gamma)$ is not unique (\cite{Ishikawa95}).

In \S \ref{Singularities of tangent varieties to surfaces.}, 
we treat the classification problem of 
singularities of tangent varieties to surfaces, exhibiting several examples and observations. 
First we observe that the tangent varieties to a generic smooth surface are not frontal. 
We characterise the surface whose tangent variety is frontal. 
In particular we show that the tangent varieties to Legendre submanifolds 
in the five dimensional standard contact projective space $P(\R^6) = \R P^5$ are frontal, if the tangent variety 
has a dense regular set. 
Recall that the singularity of tangent variety to a curve along ordinary points 
is the cuspidal edge. Therefore the singularity of tangent variety at almost any point on a curve 
is diffeomorphic to cuspidal edge, which is a generic singularity of wave front. 
We study the analogous problem for tangent varieties to Legendre surfaces. 
Then we observe that the situation becomes absolutely different. 
In fact we introduce the notion of hyperbolic and elliptic 
ordinary points on Legendre surface in  $\R P^5$ and show that 
the transverse section of the tangent variety to the surface, by a $3$-plane, has $D_4$-singularities 
(Theorem \ref{tangent to ordinary}). 

In the last section \S \ref{Tangent maps to frontal maps.}, 
we collect open problems related to several results. 

\

In this paper all manifolds and maps are assumed to be of class $C^\infty$ unless otherwise stated.



\section{Frontal maps and tangent varieties.}
\label{Frontal maps and tangent varieties.}

\bef
\label{def-frontal}
{\rm 
Let $N$ and $M$ be manifolds of dimension $n$ and $m$ respectively. 
Suppose $n \leq m$. 

A mapping $f : N^n \to M^m$ is called {\it frontal} if 
\\
{\rm (i)} the regular locus 
$$
\Reg(f) = \{ x \in N \mid f : (N, x) \to (M, f(x)) \ {\mbox{\rm is an immersion}} \}
$$ of $f$ is 
dense in $N$ and 
\\
{\rm (ii)} there exists a $C^\infty$ mapping 
$
\widetilde{f} : N \to \Gr(n, TM) = \bigcup_{y\in M} \Gr(n, T_yM)
$ 
satisfying 
$$
\widetilde{f}(x) = f_*(T_xN), \quad {\mbox{\rm for }} x \in \Reg(f). 
$$
}
\enf
Here $\Gr(n, T_yM)$ is the Grassmannian of $n$-planes in $T_yM$. 
Note that the lifting $\widetilde{f}$ is uniquely determined 
if it exists and is called the {\it Grassmannian lifting} of $f$. 

We define a subbundle ${\mathcal C} \subset T\Gr(n, TM)$ by setting, 
for $v \in T_L\Gr(n, TM), L \in T_yM$, 
$$
v \in {\mathcal C}_L \Longleftrightarrow \pi_*(v) \in L \subset T_yM. 
$$
The differential system ${\mathcal C}$ is called the {\it canonical differential system}. 
The Grassmannian lifting $\widetilde{f}$ is a ${\mathcal C}$-integral map, that is, 
$\widetilde{f}_*(TN) \subset {\mathcal C}$. We describe the canonical system in the next section 
Remark \ref{1,n+2}  
in the case $M$ is a projective space. 

If $f$ is an immersion, then $f$ is frontal. 
A wave-front hypersurface is frontal. 
The key observation for the classification of singularities of tangent varieties is 
that the tangent variety $\Tan(\gamma)$ to a curve $\gamma$ of finite type is frontal. 
the lifting Grassmannian is obtained by taking 
osculating planes to the curves (See \S \ref{Flat projective structure and the type of a curve.}). 
If $n = m$, then $f$ is frontal if the condition (i) is fulfilled, $\widetilde{f}(x)$ being $T_{f(x)}\R^m$. 

If $\widetilde{f}$ is an immersion, then the frontal mapping is called a {\it front}. 
In \cite{Ishikawa94}, we called a frontal hypersurface $(m = \ell +1)$, 
a \lq\lq front hypersurface". However 
we would like to reserve the notion \lq\lq front" 
for the case that the Grassmannian lifting is an immersion, as in the Legendre singularity theory. 
Note that frontal maps are studied also in \cite{KS}\cite{Kossowski2}\cite{SUY}. 

\bef
\label{def-tangent}
{\rm
Let $f : (\R^n, a) \to (\R^m, b), n \leq m$ be a frontal map-germ and 
$$
\widetilde{f} : (\R^n, a) \to \Gr(n, T\R^m) \cong \R^m\times \Gr(n, \R^m)
$$ be the 
Grassmannian lifting of $f$. 

A {\it tangent frame} to $f$ means a system of vector fields 
$
v_1, \dots, v_n : (\R^n, a) \to T\R^m
$ 
along $f$ such that 
$v_1(x), \dots, v_n(x)$ form a basis of $\widetilde{f}(x) \subset T_{f(x)}\R^m$. 
Then the {\it tangent map} $\Tan(f, v) : (\R^{2n}, (a, 0)) \to (\R^m, b)$ is defined by 
$$
\Tan(f, v)(s, x) := f(x) + \sum_{i=1}^n s_i v_i(x). 
$$

If we choose another tangent frame $u_1, \dots, u_n$ of $f$ and define 
$$
\Tan(f, u)(s, x) = f(x) + \sum_{i=1}^n s_i u_i(x). 
$$
Then $\Tan(f, u)$ and $\Tan(f, v)$ are right-equivalent. 
Therefore the {\it tangent variety} $\Tan(f)$ to 
a frontal map-germ is uniquely determined as a parametrised variety. 

For a frontal map-germ $f : (\R^n, 0) \to \R P^{N+1}$ in a projective space 
we define the tangent map $\Tan(f) : (\R^n, 0) \to \R P^{N+1}$ 
by taking a local projective coordinate 
$
(\R P^{N+1}, f(0)) \to (\R^{N+1}, 0)
$ 
(cf. \S \ref{Flat projective structure and the type of a curve.}). 
}
\enf

\ber
{\rm
In this paper we treat only tangent varieties, which are closely related to the {\it secant varieties}. 
The secant variety of a submanifold $S \subset \R P^n$ 
is the ruled variety obtained by taking the union of secants connecting two distinct points on $S$ and 
by taking its closure (\cite{Zak}\cite{FL}). See also Example \ref{Tangents to Veronese}. 
The secant variety is parametrised by the \lq secant map' and 
the tangent map is the \lq boundary' of secant map in some sense. 
For the singularities of secant maps, see \cite{GR}. 
}
\enr

Let $\gamma : (\R, 0) \to \R P^{N+1}$ a germ of immersion and 
$
\gamma(t) = \left( x_1(t), x_2(t), \dots, x_{N+1}(t)\right)
$
be a local representation of $\gamma$. 
Then 
$\gamma'(t)$ gives the tangent frame of $\gamma$. 
Then the tangent variety to $\gamma$ is given by 
$\Tan(\gamma) : (\R^2, 0) \to \R^{N+1}$ defined by 
$$
\Tan(\gamma)(s, t) = \gamma(t) + s\, \gamma'(t) = 
\left( x_i(t) + s\, x_i'(t)\right)_{1 \leq i \leq N+1}. 
$$
Note that $s$ is the parameter of tangent lines, while $t$ is the parameter of the original curve $\gamma$. 

If $t = 0$ is a singular point of $\gamma$, then the velocity vector $\gamma'(0) = 0$, 
and hence the above map-germ does not give the parametrisation of the tangent variety. 
However if there is $k > 0$ such that $v(t) = (1/t^k)\gamma'(t)$ is a tangent frame of $\gamma$, 
then we set 
$$
\Tan(\gamma)(s, t) = \gamma(t) + s\left(\dfrac{1}{t^k}\gamma'(t)\right) = 
\left( x_i(t) + s\left(\dfrac{1}{t^k}x_i'(t)\right)\right)_{1 \leq i \leq N+1}. 
$$
We take $k = 0$ when $\gamma$ is an immersion at $0$. 

In the above case, $\gamma$ is frontal and under a mild condition $\Tan(\gamma)$ is also frontal. 

\bet
\label{Tangent is frontal}
Let $\gamma : (\R, 0) \to \R P^{N+1}$ be a curve of finite type {\rm (}\S 
\ref{Flat projective structure and the type of a curve.}{\rm )}. 
Then $\gamma$ is frontal. Moreover the tangent map 
$\Tan(\gamma) : (\R^2, 0) \to \R P^{N+1}$ of $\gamma$ is frontal. 
\ent

Theorem \ref{Tangent is frontal} is proved in 
\S \ref{Flat projective structure and the type of a curve.}.

\ber
{\rm
Let $\gamma$ be a curve of finite type. Then 
it is natural to ask what $\Tan(\Tan(\gamma))$ is, because 
$\Tan(\gamma)$ is frontal. For a curve $\gamma$ in $\R P^{N+1}, N \geq 2$, 
the tangent plane to $\Tan(\gamma)$ along each ruling (tangent line) 
is constant, that is the osculating $2$-plane. Therefore $\Tan(\Tan(\gamma))$ 
is a $3$-fold, not a $4$-fold, ruled by osculating $2$-planes of the original curve $\gamma$ (\cite{Ishikawa95}). 
}
\enr

\

We classify the map-germ $\Tan(\gamma)$ by local right-left diffeomorphism equivalence. 
Two map-germs $f : (N, a) \to (M, b)$ and 
$f' : (N', a') \to (M', b')$ are called {\it diffeomorphic} or {\it right-left equivalent} if 
there exist diffeomorphism-germs $\sigma : (N, a) \to (N', a')$ and 
$\tau : (M, b) \to (M', b')$ such that 
$f'\circ\sigma = \tau\circ f$. 

In the followings, $I$ is an open interval. 

\bet {\ }
\label{tangent to generic curve}
{\rm (1)} {\rm (\cite{Cleave})} 
For a generic curve $\gamma : I \to \R P^3$ in $C^\infty$-topology, 
the curve $\gamma$ is of finite type at each point in $I$ and 
the tangent variety $\Tan(\gamma)$ to $\gamma$ at each point in $I$ 
is locally diffeomorphic to the cuspidal edge or to 
the folded umbrella {\rm (}cuspidal cross cap{\rm)}. 

{\rm (2)} 
Let $N+1 \geq 4$. For a generic curve $\gamma : I \to \R P^{N+1}$ in $C^\infty$-topology, 
the curve $\gamma$ is of finite type at each point in $I$ and 
the tangent variety $\Tan(\gamma)$ to $\gamma$ at each point of $I$ 
is locally diffeomorphic to the cuspidal edge. 
\ent

The genericity means the existence of an open dense subset ${\mathcal O} \subset C^\infty(I, \R P^{N+1})$ 
such that any $\gamma \in {\mathcal O}$ satisfies the consequence. 

The {\it cuspidal edge} is parametrised by the map-germ $(\R^2, 0) \to (\R^{N+1}, 0), (N+1 \geq 3)$ defined by 
$$
(u, x) \mapsto (u, \ x^2, \ x^3, \ 0, \dots, \ 0). 
$$
Note that it is diffeomorphic (right-left equivalent) to the germ 
$$
(t, s) \mapsto (t+s, \ t^2 + 2s t, \ t^3 + 3s t^2, \ \dots, \ t^{N+1} + (N+1)s t^N), 
$$
and also to 
$$
(t, s) \mapsto (t+s, \ t^2 + 2s t, \ t^3 + 3s t^2, \ 0, \dots, \ 0), 
$$
A {\it folded umbrella} is parametrised by the germ $(\R^2, 0) \to (\R^3, 0)$ defined by 
$$
(t, s) \mapsto (t + s, t^2 + 2s t, t^4 + 4s t^3), 
$$
which is diffeomorphic to 
$$
(u, x) \mapsto (u, \ x^2 + ux, \ \dfrac{1}{2}x^4 + \dfrac{1}{3}ux^3). 
$$
A folded umbrella is often called a cuspidal cross cap. 

\begin{figure}[htbp]
\begin{center} 
\includegraphics[width=7truecm, height=2truecm, clip, 
bb=20.9189 738.12 311.453 818.127]{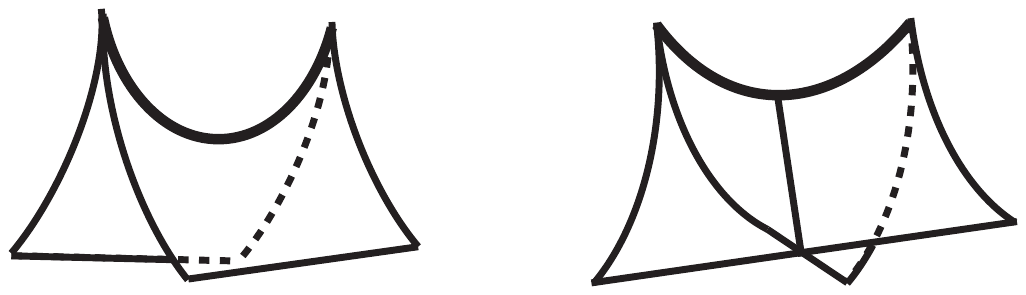}  
\caption{cuspidal edge and folded umbrella.}
\label{3-space}
\end{center}
\end{figure}

Theorem \ref{tangent to generic curve} is proved in \S \ref{Normal forms of tangent surfaces.}.

\bee
\label{umbilical bracelet}
{\rm (umbilical bracelet)
Let 
$$
V^{N+2} = \{ a_0 x^{N+1} + a_1 x^Ny + \cdots + a_N xy^N + a_{N+1} y^{N+1} \} \cong \R^{N+2}
$$ 
be the space of homogeneous polynomials of degree $N+1$ in two variables $x, y$. 
The polynomials with zeros of multiplicity $N+1$ form a curve $C$ in $P(V) \cong \R P^{N+1}$. 
The tangent variety $\Tan(C)$ to $C$ 
coincides with the set of polynomials with zeros of multiplicity $\geq N$. 
The surface $\Tan(C)$  has cuspidal edge singularities along $C$. 
In particular in the case $N+1 = 3$, 
the tangent variety $\Tan(C)$ to $C$ is called the {\it umbilical bracelet}(\cite{Porteous}\cite{FW}).  
If $N+1 \geq 4$, $\Tan(\Tan(C)) \subset P(V^{N+2})$ 
coincides with of polynomials with with zeros of multiplicity $\geq N-1$. 
}
\ene

\ber
{\rm
The tangent surface to a curve is obtained as a union of strata of envelope generated by the 
dual curve to the original curve. 
The generating family associated to the dual curve is determined, up to parametrised ${\mathcal K}$-equivalence 
in several cases. 
We recall the notion of types of curves in a projective space 
in \S \ref{Flat projective structure and the type of a curve.}. 
If the type ${\mathbf A} = (a_1, \dots, a_{N+1})$ of a curve in $\R P^{N+1}$ 
is one of followings 
$$
\begin{array}{lcl}
{\rm (I)}_{N,r} & : & (1, 2, \dots, N, N+r), \ (r = 0, 1, 2, \dots), 
\\
{\rm (II)}_{N,i} & : & (1, 2, \dots, i, i+2, \dots, N+1, N+2), \ (0 \leq i \leq N-1), 
\\
{\rm (III)}_{N} & : & (3, 4, \dots, N+2, N+3), 
\end{array}
$$
then the generating family is determined by the type of the curve \cite{Ishikawa93}. 
In each case, a normal form of the tangent variety can be obtained from 
the generating family 
$$
F(t, x) = t^{a_{N+1}} + x_1t^{a_{N+1}-a_1} + x_2 t^{a_{N+1}-a_2} + \cdots + x_N t^{a_{N+1}-a_N} + x_{N+1} = 0, 
$$
by solving 
$$
F = 0, \quad \dfrac{\pa F}{\pa t} = 0, \quad \dots, \quad \dfrac{\pa^{N-1} F}{\pa t^{N-1}} = 0, 
$$
deleting the divisor $\{t = 0\}$ if necessary. 
For example, for the type ${\rm (II)}_{3,2} : (1, 2, 4, 5)$, 
we have generating family 
$$
F(t, x) = t^5 + x_1 t^4 + x_2 t^3 + x_3 t + x_4. 
$$
Then the tangent variety is obtained by solving 
$$
\left\{ 
\begin{array}{l}
t^5 + x_1t^4 + x_2 t^3 + x_3 t + x_4 = 0, 
\\
5t^4 + 4x_1t^3 + 3x_2 t^2 + x_3 = 0, 
\\
20t^3 + 12x_1 t^2 + 6x_2 t = 0. 
\end{array}
\right.
$$
In fact, from these equations, we get a map-germ $(\R^2, 0) \to (\R^4, 0)$ by 
$$
x_2 = -\dfrac{10}{3}t^2 - 2x_1t, \quad 
x_3 = 5t^4 + 2x_1t^3, \quad 
x_4 = -\dfrac{8}{3}t^5 - x_1t^4, 
$$
which is diffeomorphic to the open folded umbrella (see Theorems \ref{osculating-frame-geq4}, \ref{1,2,4,5}). 
}
\enr

\section{Differential systems on flag manifolds.}
\label{Differential systems on flag manifolds.}

First we recall the flag manifolds and the canonical differential systems on flag manifolds. 
For the generality on differential systems, see \cite{IL}. 

Let $V$ be a vector space of dimension $n$ and 
$0 < n_1 < n_2 < \cdots < n_\ell < n$. 
Then we define the flag manifold 
$$
{\mathcal F} = {\mathcal F}_{n_1, n_2, \dots, n_\ell}(V) := 
\left\{ V_{n_1} \subset V_{n_2} \subset \cdots \subset V_{n_\ell} \subset V \mid 
\dim(V_{n_j}) = n_j, (1 \leq j \leq \ell) \right\}. 
$$
Note that 
$$
\dim({\mathcal F}) = n_1(n - n_1) + (n_2 - n_1)(n - n_2) + \cdots + (n_\ell - n_{\ell-1})(n - n_\ell).
$$

Denote by $\pi_i : {\mathcal F}_{n_1, n_2, \dots, n_\ell}(V) \to 
{\Gr}(n_i, V)$ the canonical projection to the $i$-th member of the flag. 
The canonical differential system ${\mathcal C} = {\mathcal C}_{n_1, n_2, \dots, n_\ell} \subset 
T{\mathcal F}$ is defined by, for $v \in T_{\mathbf V}{\mathcal F}, {\mathbf V} \in {\mathcal F}$, 
$$
v \in {\mathcal C}_{\mathbf V} \Longleftrightarrow  
{\pi_i}_*(v) \in T{\Gr}(n_i, V_{n_{i+1}}) (\subset T{\Gr}(n_i, V)), 
(1 \leq i \leq \ell-1). 
$$
Then ${\mathcal C}$ is a bracket-generating (completely non-integrable) subbundle of $T{\mathcal F}$ with 
$$
\rank({\mathcal C}) = n_1(n_2 - n_1) + (n_2 - n_1)(n_3 - n_2) + \cdots + (n_\ell - n_{\ell-1})(n - n_\ell). 
$$

A $C^\infty$ curve $\Gamma : I \to {\mathcal F}$ from an open interval $I$ is called
a ${\mathcal C}$-{\it integral curve} if 
$\Gamma'(t) \in {\mathcal C}_{\Gamma(t)}$ for any $t \in I$. 
A ${\mathcal C}$-integral curve can be phrased as a $C^\infty$-family 
$$
c(t) = (V_{n_1}(t), V_{n_2}(t), \dots, V_{n_\ell}(t))
$$ 
of flags in ${\mathcal F}$ such that 
each $V_{n_i}(t)$ moves along $V_{n_{i+1}}(t)$ at every moment infinitesimally. 


\

Let $V$ be an $(N+2)$-dimensional vector space. 
For the study of tangent varieties to curves, it is natural to regard the following flag manifolds 
$$
{\mathcal F}_{1,2} = {\mathcal F}_{1,2}(V) := 
\{ V_1 \subset V_2 \subset  V \mid \dim(V_i) = i \}, 
$$
and 
$$
{\mathcal F}_{1,2,3} = {\mathcal F}_{1,2,3}(V) := 
\{ V_1 \subset V_2 \subset V_3 \subset  V \mid \dim(V_i) = i \}. 
$$
The canonical systems 
${\mathcal T} = {\mathcal C}_{1,2}$ and ${\mathcal N} = {\mathcal C}_{1,2,3}$ are defined as follows: 
For $(V_1, V_2) \in {\mathcal F}_{1,2}$,  
$$
v \in {\mathcal T}_{(V_1, V_2)} \Longleftrightarrow  
{\pi_1}_*(v) \in TP(V_2) (\subset TP(V)). 
$$ 
For $(V_1, V_2, V_3) \in {\mathcal F}_{1,2,3}$, 
$$
w \in {\mathcal N}_{(V_1, V_2, V_3)} \Longleftrightarrow  
{\pi_1}_*(w) \in TP(V_2) (\subset TP(V)), \ {\pi_2}_*(w) \in T\Gr(2, V_3) (\subset T\Gr(2, V)). 
$$

Then we have 

\bep
Let $\gamma : (\R, 0) \to P(V^{N+2}) \cong \R P^{N+1}$ be a $C^\infty$ curve. Suppose 
$\Reg(\gamma)$ is dense in $(\R, 0)$. 
Then $\gamma$ is frontal if and only if $\gamma = \pi_1\circ c$ for some 
${\mathcal C}_{1,2}$-integral curve $c : (\R, 0) \to {\mathcal F}_{1,2}(V)$. 
\enp

In fact $c$ gives a tangent frame of $\gamma$. In this case, $\gamma$ is called {\it tangent-framed}. 

\bep
\label{1,2,3--frontal}
Let $\gamma : (\R, 0) \to P(V^{N+2}) \cong \R P^{N+1}$ be a frontal curve. Suppose 
$\Reg(\Tan(\gamma))$ is dense in $(\R^2, 0)$. 
Then $\Tan(\gamma)$ is frontal if and only if $\gamma = \pi_1\circ \kappa$ for some 
${\mathcal C}_{1,2,3}$-integral curve $\kappa : (\R, 0) \to {\mathcal F}_{1,2,3}(V)$. 
\enp

In fact, if $\Tan(\gamma)$ is frontal, then $V_1(t) = \gamma(t)$, 
the tangent (projective) line $V_2(t)$ to $\gamma$ at $t$ and the tangent (projective) 
plane to $\Tan(\gamma)$ 
at $(t, 0)$ form a ${\mathcal C}_{1,2,3}$-integral lifting of $\gamma$. Conversely 
if $c(t) = (V_1(t), V_2(t), V_3(t))$ is a ${\mathcal C}_{1,2,3}$-integral curve, then 
$\Tan(\gamma)$ has the constant tangent plane $V_3(t)$ along each ruling, and $(t, s) 
\mapsto T_{V_1(t)}P(V_3)(t)$ gives the Grassmannian lifting of $\Tan(\gamma)$. 

The projection of a ${\mathcal C}_{1,2,3}$-integral curve 
is called a {\it tangent-principal-nomal-framed} curve. 

%

\bet 
\label{tangent-framed classification}
{\rm (1)}
Let $N+1=3$. 
For a generic ${\mathcal C}_{1,2}$-integral curve $c : I \to {\mathcal F}_{1,2}(V^{4})$ in $C^\infty$-topology,  
the tangent variety $\Tan(\gamma)$ 
to the tangent-framed curve $\gamma = \pi_1\circ c : I \to P(V^4) = \R P^3$ at each point 
is locally diffeomorphic to the cuspidal edge, the folded umbrella or the swallowtail. 

{\rm (2)}
Let $N+1 \geq 4$. 
For a generic ${\mathcal C}_{1,2}$-integral curve $c : I \to {\mathcal F}_{1,2}(V^{N+2})$ in $C^\infty$-topology,  
the tangent variety $\Tan(\gamma)$ 
to the tangent-framed curve $\gamma = \pi_1\circ c : I \to P(V) = \R P^{N+1}$ at each point 
is locally diffeomorphic to the cuspidal edge or the open swallowtail. 
\ent

The {\it swallowtail} $(\R^2, 0) \to (\R^3, 0)$ is given by 
$$
(t, s) \mapsto (t^2 + 2s, \ t^3 + 3st, \ t^4 + 4st^2), 
$$
which is diffeomorphic to 
$$
(u, x) \mapsto (u, \ x^3 + ux, \ \frac{3}{4}x^4 + \frac{1}{2}ux^2).
$$
The 
{\it open swallowtail} $(\R^2, 0) \to (\R^{N+1}, 0), N+1 \geq 4$ is given by 
$$
(t, s) \mapsto (t^2 + 2s, \ t^3 + 3st, \ t^4 + 4st^2, \ t^5 + 5st^3, \ 0, \ \dots, \ 0), 
$$
which is diffeomorphic to 
$$
(u, x) \mapsto (u, \ x^3 + ux, \ \frac{3}{4}x^4 + \frac{1}{2}ux^2, \ \frac{3}{5}x^5 + \frac{1}{3}ux^3, \ 0, \ \dots, \ 0). 
$$

\bet 
\label{tangent-principal-normal-framed classification}
{\rm (1)}
Let $N+1=3$. 
For a generic ${\mathcal C}_{1,2,3}$-integral curve $\kappa : I \to {\mathcal F}_{1,2,3}(V^{4})$ in $C^\infty$-topology,  
the tangent variety $\Tan(\gamma)$ 
to the tangent-principal-normal-framed curve $\gamma = \pi_1\circ \kappa : I \to P(V^4) = \R P^3$ at each point 
is locally diffeomorphic to the cuspidal edge, the folded umbrella, the Mond surface or the swallowtail. 

{\rm (2)}
Let $N+1 \geq 4$. 
For a generic ${\mathcal C}_{1,2,3}$-integral curve $\kappa : I \to {\mathcal F}_{1,2,3}(V^{N+2})$ 
in $C^\infty$-topology,  
the tangent variety $\Tan(\gamma)$ 
to the tangent-principal-normal-framed 
curve $\gamma = \pi_1\circ \kappa : I \to P(V) = \R P^{N+1}$ at each point 
is locally diffeomorphic to the cuspidal edge, the open Mond surface or the open swallowtail. 
\ent

The {\it Mond surface} $(\R^2, 0) \to (\R^3, 0)$ is given by 
$$
(t+s, \ t^3 + 3st^2, \ t^4 + 4st^3), 
$$
which is diffeomorphic to 
$$
(u, x) \mapsto (u, \ x^3 + ux^2, \ \frac{3}{4}x^4 + \frac{2}{3}ux^3)). 
$$
The Mond surface is called also {\it cuspidal beaks} (\cite{IS})
or {\it cuspidal beak to beak} (\lq bec \`{a} bec'). 

The {\it open Mond surface} $(\R^2, 0) \to (\R^{N+1}, 0), N+1 \geq 4$ is given by 
$$
(t, s) \mapsto (t + s, \ t^3 + 3st^2, \ t^4 + 4st^3, \ t^5 + 5st^4, \ 0, \ \dots, \ 0), 
$$
$$
(u, x) \mapsto (u, \ x^3 + ux^2, \ \frac{3}{4}x^4 + \frac{2}{3}ux^3, \ \frac{3}{5}x^5 + \frac{1}{2}ux^4, \ 0, \ \dots, \ 0). 
$$

\

Now we recall on osculating-framed curves (cf. \cite{Ishikawa12}). 
Let $V$ be an $(N+2)$-dimensional real vector space. 
Consider the complete flag manifold: 
$$
{\mathcal F} = {\mathcal F}_{1,2,\dots,N+1}(V) := \{ V_1 \subset V_2 \subset \cdots V_{N+1} \subset V \mid 
\dim(V_i) = i, 1 \leq i \leq N+1\}. 
$$
Then $\dim {\mathcal F} = \frac{(N+1)(N+2)}{2}$. 
We denote by $\pi_i : {\mathcal F} \to \Gr(i, V)$ the canonical projection 
$$
\pi_i(V_1, V_2, \dots, V_{N+1}) = V_i. 
$$

The canonical system ${\mathcal C} = {\mathcal C}_{1,2,\dots,N+1} 
\subset T{\mathcal F}$ 
is defined by  
$$
v \in {\mathcal C}_{(V_1, \dots, V_{N+1})} \Longleftrightarrow  
{\pi_i}_*(v) \in T{\Gr}(i, V_{i+1}) (\subset T{\Gr}(i, V)), 
(1 \leq i \leq N). 
$$



For a $C^\infty$ curve $\gamma : I \to P(V) = \R P^{N+1}$, 
if we consider Frenet-Serret frame, or the osculating projective moving frame, 
$
\Gamma = (e_0(t), e_1(t), \dots, e_{N+1}(t)) : 
I \to {\mbox{\rm GL}}(\R^{N+2}) = {\mbox{\rm GL}}(N+2, \R), 
$
$\gamma(t) = [e_0(t)]$, 
then, setting 
$
V_i(t) := \langle e_0(t), e_1(t), \dots, e_{i-1}(t) \rangle_{\R}, (1 \leq i \leq N+1), 
$
we have a ${\mathcal C}$-integral lifting $\widetilde{\gamma} : I \to {\mathcal F}$ 
of $\gamma$ for the projection $\pi_1 : {\mathcal F} \to P(V)$, 
by 
$\widetilde{\gamma}(t) =  (V_1(t), V_2(t), \dots, V_{N+1}(t))$. 
%
In this case, $\gamma$ is called {\it osculating-framed}. 
Note that 
the framing of an osculating-framed curve is uniquely determined if an orientation of the curve 
and a metric on $P(V)$ are given. 

\begin{figure}[htbp]
\begin{center} 
\includegraphics[width=12truecm, height=2.5truecm, clip, 
bb=20.9189 725.354 569.037 822.88]{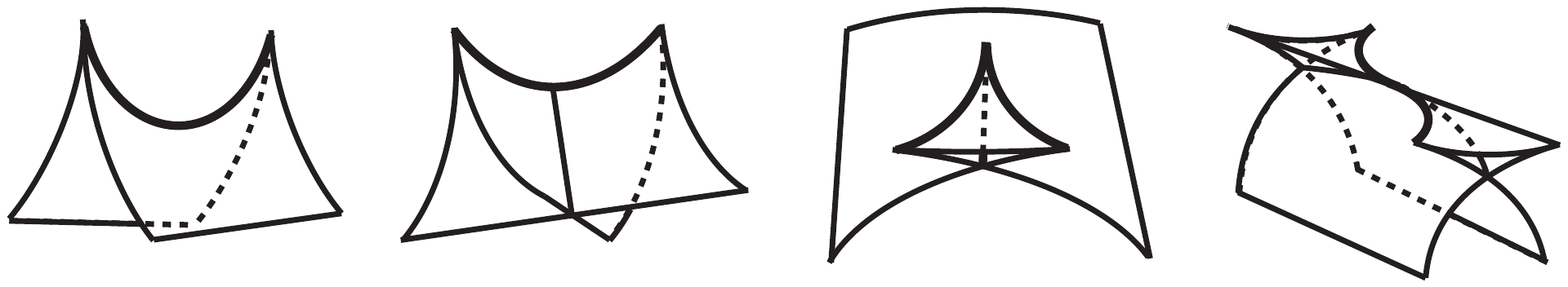}  
\caption{cuspidal edge, folded umbrella, swallowtail and Mond surface in $\R^3$.}
\label{3-space-frame}
\end{center}
\end{figure}

\bet 
{\rm (\cite{Ishikawa12})} Let $N+1 = 3$. 
For a generic ${\mathcal C}$-integral curve $c : I \to {\mathcal F}(V^4)$ in $C^\infty$-topology,  
the tangent variety $\Tan(\gamma)$ 
to the osculating-framed curve $\gamma = \pi_1\circ c : I \to P(V^4) = \R P^3$ at each point of $I$ 
is locally diffeomorphic to the cuspidal edge, the folded umbrella, 
the swallowtail or to the Mond surface {\rm (Figure \ref{3-space-frame})}.
\ent

In this paper we treat higher codimensional cases, and we show the following

\begin{figure}[htbp]
\begin{center}
 \includegraphics[width=12truecm, height=2.5truecm, clip, bb=5.9165 696.727 582.302 807.316]{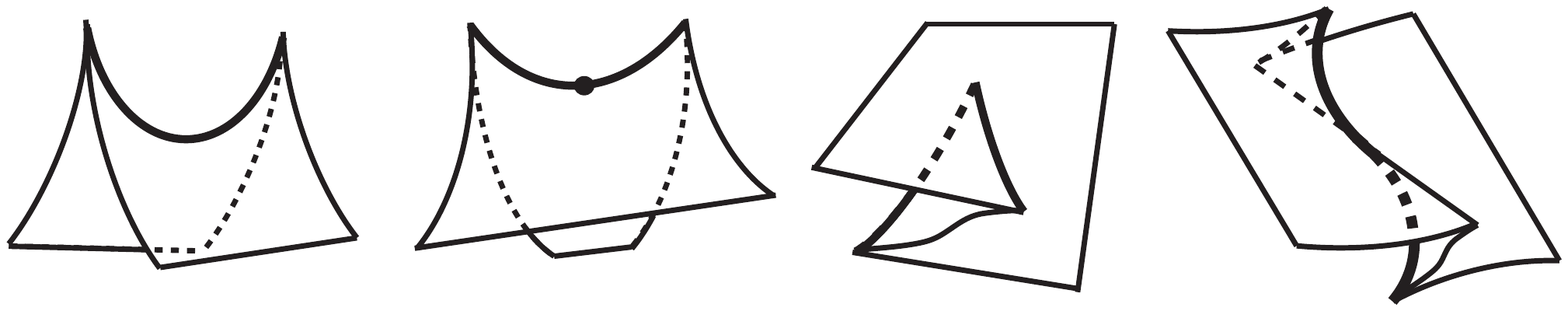}  
\caption{cuspidal edge, open folded umbrella, open swallowtail and open Mond surface in $\R^4$.}
\label{higher-framed}
\end{center}
\end{figure}

\bet
\label{osculating-frame-geq4}
Let $N+1 \geq 4$. 
For a generic ${\mathcal C}$-integral curve $c : I \to {\mathcal F}$ in $C^\infty$-topology,  
the tangent variety to the osculating-framed curve $\gamma = \pi_1\circ c : I \to P(V^{N+2}) = \R P^{N+1}$ 
at each point is locally diffeomorphic to the cuspidal edge, the open folded umbrella {\rm (}cuspidal non-cross cap{\rm )}, 
the open swallowtail or to the open Mond surface
{\rm (Figure \ref{higher-framed})}.
\ent

The {\it open folded umbrella} $(\R^2, 0) \to (\R^{N+1}, 0), N \geq 3$ is given by 
$$
(t, s) \mapsto (t + s, \ t^2 + 2st, \ t^4 + 4st^3, \ t^5 + 5st^4, \ 0, \ \dots, \ 0), 
$$
which is diffeomorphic to 
$$
(u, x) \mapsto (u, \ x^2 + ux, \ \frac{1}{2}x^4 + \frac{1}{3}ux^3, \ \frac{2}{5}x^5 + \frac{1}{4}ux^4, \ 0, \ \dots, \ 0). 
$$

Our main results, Theorems \ref{tangent-framed classification}, \ref{tangent-principal-normal-framed classification}, 
\ref{osculating-frame-geq4} are proved in \S \ref{Normal forms of tangent surfaces.}. 

\

Lastly in this section, we describe the 
canonical system ${\mathcal C} = {\mathcal C}_{1,2,\dots,k+1}$ on ${\mathcal F}_{1,2,\dots,k+1}(V^{N+2})$. 
Let ${\mathbf V}_1 = (V_{11}, V_{21}, \dots, V_{k+1,1}) \in {\mathcal F}_{1,2,\dots,k+1}(V^{N+2})$. 
Fix a flag $V^{N+2} \supset W_{N+1} \supset W_N \supset W_{N-k+1}$ such that 
$W_{N-i+1}\cap V_{i+1\ 1} = \{ 0\}$, $i = 0, 1, \dots, k$. Take the open neighbourhood $U$ of ${\mathbf V}_1$ 
defined by 
$$
U := \left\{ (V_1, V_2, \dots, V_{k+1}) \in {\mathcal F}_{1,2,\dots,k+1}(V^{N+2}) \ \ \left\vert 
\ \ 
W_{N-i+1}\cap V_{i+1} = \{ 0\}, \ i = 0, 1, \dots, k\right\}.\right.
$$
Take non-zero vectors $e_0 \in V_{11}, e_1 \in V_{21}\cap W_{N+1}, e_2 \in V_{31}\cap W_N, 
\dots, e_{k} \in V_{k+1\, 1}\cap W_{N-k+2}$. Adding a basis $(e_{k+1}, \dots, e_{N+1})$ of $W_{N-k+1}$, we 
get a basis $(e_0, e_1, e_2, \dots, e_k, e_{k+1}, \dots, e_{N+1})$ of $V$. 
Then, for each ${\mathbf V} = (V_1, V_2, \dots, V_{k+1})$, $V_{i}$ has a basis $v_0, v_1, \dots, v_{i-1}$ 
(a \lq moving frame') of the form 
$$
v_i = e_i + \sum_{j=i+1}^{N+1} x_j^{\ i}e_j, \quad 0 \leq i \leq k. 
$$
Then the condition that a curve in ${\mathcal F}_{1,2,\dots,k+1}(V^{N+2})$ ${\mathcal C}$-integral 
is equivalent to that 
the components of the curve satisfies the conditions  
$$
(v_{i-1})' = \sum_{j=i}^{N+1} (x_j^{\ i-1})'e_j \in \langle v_0, v_1, \dots, v_{i}\rangle_{{\mathcal E}_1}, \quad 
1 \leq i \leq k. 
$$
Thus we see that the differential system ${\mathcal C} = {\mathcal C}_{1,2,\dots,k+1}$ 
is defined by 
$$
dx_j^{\ i-1} - x_j^{\ i} dx_{i}^{\ i-1} = 0, \quad (1 \leq i \leq k, i+1 \leq j \leq N+1), 
$$
for the system of local coordinates $\left( x_j^{\ i}\right)_{0 \leq i \leq k, i+1 \leq j \leq N+1}$ of 
${\mathcal F}_{1,2,\dots,k+1}(V^{N+2})$. 

\

\ber
\label{1,n+2}
{\rm 
For a $(N+2)$-dimensional vector space $V$, 
the Grassmannian bundle 
$\Gr(n, TP(V^{N+2}))$ over $P(V^{N+2})$ is 
identified with the flag manifold ${\mathcal F}_{1,n+1}(V^{N+2})$, 
$$
{\mathcal F}_{1,n+1}(V^{N+2}) = \{ V_1 \subset V_{n+1} \subset V^{N+2} \mid \dim(V_1) = 1, \dim(V_{n+1}) = n+1 \}. 
$$
Remark that 
the Grassmannian liftings of frontal maps $N^n \to P(V^{N+2})$ are ${\mathcal C}$-integral 
of the canonical system ${\mathcal C} = {\mathcal C}_{1,n+1}$. 

The canonical system ${\mathcal C}_{1,n+1}$ on ${\mathcal F}_{1,n+1}(V^{N+2})$ is locally given by 
$$
dx_{i+1}^{\ 0} - \sum_{j=1}^n x_{i+1}^{\ j}dx_j^{0} = 0, \quad (n \leq i \leq N), 
$$
for a system of local coordinates $x_{i+1}^{\ 0}, (0 \leq i \leq N), x_{i+1}^{\ j}, (1 \leq j \leq n, 
n \leq i \leq N)$. The projection $\pi_1 : \Gr(n, TP(V^{N+2})) \to P(V^{N+2})$ 
is represented by $(x_1^{\ 0}, \dots, x_{N+1}^{\ 0})$. 
If we write $x_i = x_{i}^{\ 0}\  (1 \leq i \leq n), y_k = x_{n+k}^{\ 0}\ (1 \leq k \leq N-n+1)$ and 
$p_k^{\ i} = x_{n+k}^{\ i}\ (1 \leq k \leq N-n+1, 1 \leq i \leq n)$, then we have
$$
dy_k - \sum_{i=1}^n p_k^{\ i} dx_i = 0, \quad 1 \leq k \leq N-n+1. 
$$
Therefore the condition that a map $F : L^n \to \Gr(n, TP(V^{N+2}))$ is ${\mathcal C}$-integral 
is expressed by 
$$
d(y_k\circ F) - \sum_{i=1}^n (p_k^{\ i}\circ F)\, d(x_i\circ F) = 0, \quad 1 \leq k \leq N-n+1. 
$$
}
\enr

\section{Type of a curve in a space with flat projective structure.}
\label{Flat projective structure and the type of a curve.}

Let $M$ be an $m$-dimensional $C^\infty$ manifold. 
A {\it flat projective structure} on $M$ is given by an atlas 
$\{(U_\alpha, \varphi_\alpha)\}$ where $M = \bigcup_{\alpha} U_\alpha$, 
$\varphi_\alpha : U_\alpha \to \varphi_\alpha(U_\alpha) \subset \R^m$, 
and transition functions 
$\varphi_{\beta}\circ \varphi_{\alpha}^{-1} : \varphi_\alpha(U_{\alpha}\cap U_{\beta}) 
\to \varphi_{\beta}(U_{\alpha}\cap U_{\beta})$ are fractional linear with a common denominator. 
Then an admissible chart is called a {\it system of projective local coordinates}.  
The projective space $P(V^{m+1})$ for a vector space $V^{m+1}$ has the canonical 
flat projective structure. 
Also Grassmannians and Lagrange Grassmannians have flat projective structures (cf. 
\cite{IMT}). 

Let $\gamma : I \to M$ be a $C^\infty$-curve in a manifold $M$ with a flat projective structure. 
Take a system of projective local coordinates $(x_1, x_2, \dots, x_m)$ 
centred at $\gamma(t_0)$ and the local affine representation $(\R, t_0) \to (\R^m, 0)$, 
$$
\gamma(t) = {}^T(x_1(t), x_2(t), \dots, x_m(t))
$$ of $\gamma$. 
Consider the $(m\times k)$-matrix
$$
W_k(t) := \left( \gamma'(t_0), \gamma''(t_0), \cdots, \gamma^{(k)}(t_0)\right) 
$$
for any integer $k \geq 1$ and $k = \infty$. 
Note that the rank of $W_k(t_0)$ is independent of the choice on representations for $\gamma$. 

\bef
\label{type-def}
{\rm 
We call $\gamma$ {\it of finite type} at $t = t_0 \in I$ if 
the $(m\times\infty)$-matrix 
$$
W_{\infty}(t_0) = \left( \gamma'(t_0), \ \gamma''(t_0), \ \cdots, \ \gamma^{(k)}(t_0), \ \cdots\cdots \ \right)
$$
is of rank $m$. 
Define, for $1 \leq i \leq m$, 
$
a_i := \min\left\{ k \mid W_k(t_0) = i\right\}.  
$
Then we have a sequence of natural numbers $1 \leq a_1 < a_2 < \dots < a_m$, 
and we call $\gamma$ {\it of type} $(a_1, a_2, \dots, a_m)$ at $t = t_0 \in I$.  

If $(a_1, a_2, \dots, a_m) = (1, 2, \dots, m)$, then 
$t = t_0$ is called an {\it ordinary point} of $\gamma$. 
}
\enf

It is easy to see 

\bel
\label{affine coord}
A curve-germ $\gamma : (\R, 0) \to M$ in a manifold $M$ with a flat projective structure, 
is of type $(a_1, a_2, \dots, a_m)$ at $0$ if and only if 
there exists a system of projective local coordinates $(x_1, x_2, \dots, x_m)$ 
centred at $\gamma(0)$ such that 
$$
x_1(t) = t^{a_1} + o(t^{a_1}), \ 
x_2(t) = t^{a_2} + o(t^{a_2}),\ \dots, \
x_m(t) = t^{a_m} + o(t^{a_m}). 
$$
\enl

\bel
Let $\gamma : (\R, 0) \to P(\R^{N+2}) = \R P^{N+1}$ be a curve and 
$\widetilde{\gamma} : (\R, 0) \to \R^{N+2} \setminus \{ 0\}$ be a lifting of $\gamma$. 
Set $\widetilde{W}_r(t) = (\widetilde{\gamma}(t), \widetilde{\gamma}'(t), \cdots,  \widetilde{\gamma}^{(r)}(t) )$. 
Then $\gamma$ is of type $\AAA = (a_1, a_2, \dots, a_{N+1})$ if and only if 
$a_i = \min \{ r \mid \widetilde{W}_r(t_0) = i+1 \}, 1 \leq i \leq N+1$.  
\enl

Moreover we see

\bel
\label{integral lifting}
Let 
$\gamma : (\R, 0) \to P(\R^{N+2}) = \R P^{N+1}$ a curve of finite type. 
There is unique ${\mathcal C}_{1,2,3,\dots,N+1}$-integral $C^\infty$ lifting $\Gamma : (\R, 0) \to 
{\mathcal F}_{1,2,3,\dots,N+1}(\R^{N+2})$ of $\gamma$. Moreover 
by the projection of $\Gamma$, we have 
${\mathcal C}_{1,2,3}$-integral lifting $\kappa : (\R, 0) \to 
{\mathcal F}_{1,2,3}(\R^{N+1})$ and ${\mathcal C}_{1,2}$-integral lifting $c : (\R, 0) \to 
{\mathcal F}_{1,2}(\R^{N+1})$ of $\gamma$. 
\enl

\Proof
The first half is proved in \cite{Ishikawa12} (Lemma 6.1). 
We take 
the lifting $\widetilde{\gamma} : (\R, 0) \to \R^{N+2} \setminus \{ 0\}$ defined by 
$$
\widetilde{\gamma}(t) = {}^T\left(1, \ t^{a_1} + o(t^{a_1}), \ t^{a_2} + o(t^{a_2}),\ \dots, \ t^{a_{N+1}} + o(t^{a_{N+1}})\right). 
$$
of $\gamma$. Consider the $(N+2)\times(N+2)$-matrix 
$$
A(t) = \left( 
\widetilde{\gamma}(t), \dfrac{1}{a_1!}\widetilde{\gamma}^{(a_1)}(t), \cdots, 
\dfrac{1}{a_{N+1}!}\widetilde{\gamma}^{(a_{N+1})}(t)
\right). 
$$
Let $V_i(t)$ be the linear subspace of $\R^{N+2}$ generated by the first $i$-columns of $A(t)$. 
Then $\Gamma : (\R, 0) \to {\mathcal F}_{1,2,3,\dots,N+1}(\R^{N+2})$ 
is uniquely determined by $\Gamma(t) = (V_1(t), V_2(t), \dots, V_{N+1}(t))$. 
The lower triangle components of $A(t)$ give the local representation of $\Gamma$, 
therefore $\Gamma$ is $C^\infty$. Moreover 
$\kappa(t) = (V_1(t), V_2(t), V_3(t))$ and $c(t) = (V_1(t), V_2(t))$. 
\QED

\

\noindent 
{\it Proof of Theorem \ref{Tangent is frontal} :} 
Theorem \ref{Tangent is frontal} follows from Lemma \ref{integral lifting} and 
Proposition \ref{1,2,3--frontal}. 
Here we give concretely the Grassmannian lifting of $\Tan(\gamma)$ in term of Wronskian. 

\bel
\label{Wronskian}
Let $\gamma : (\R, 0) \to \R P^{N+1}$ 
be a curve-germ of type $(a_1, a_2, \dots, a_{N+1})$ 
and 
$$
\gamma(t) = \left( x_1(t), x_2(t), \dots, x_{N+1}(t)\right)
$$
be a local affine representation of $\gamma$. 
Then the tangent variety to $\gamma$ is parametrised by 
$$
f(s, t) = \Tan(\gamma)(s, t) := \gamma(t) + s\, \frac{1}{\alpha(t)}\gamma'(t) = 
\left( x_i(t) + s\, \frac{1}{\alpha(t)}x_i'(t)\right)_{1 \leq i \leq N+1}, 
$$
where $\alpha(t) = t^{a_1 - 1}$. 
We set $f_i(s, t) = x_i(t) + \dfrac{s}{\alpha(t)}x_i'(t)$. 
Then we have 
$$
\frac{W_{i2}}{W_{12}} \, df_1 + \frac{W_{1i}}{W_{12}} \, df_2. 
$$
Here $$
W_{ij}(t) = 
\left\vert
\begin{array}{cc}
x_i'(t) & x_j'(t) \\
x_i''(t) & x_j''(t) 
\end{array}
\right\vert.
$$
\enl

\Proof 
We have 
$$
df_i(s, t) = 
\frac{x_i'(t)}{\alpha(t)} \,ds  + \left( x_i'(t) + s\, \left(\frac{x_i'(t)}{\alpha(t)}\right)'\,\right)dt. 
$$
Then we have 
$$
\left\vert
\begin{array}{cc}
x_1' & x_2' \\
x_1' + s\left(\frac{x_1'}{\alpha}\right)' & x_2' + s\left(\frac{x_2'}{\alpha}\right)'
\end{array}
\right\vert
= 
\dfrac{s}{\alpha}\,W_{12}. 
$$
Therefore we have, for $3 \leq i \leq n+1$,  
$$
\begin{array}{rcl}
df_i & = & 
{\displaystyle 
\frac{\alpha}{sW_{12}}
\left( df_1, \ df_2\right)
\left(
\begin{array}{cc}
x_2' + s\left(\frac{x_2'}{\alpha}\right)' & -x_2' \\
-x_1' - s\left(\frac{x_1'}{\alpha}\right)' & x_1'
\end{array}
\right)
\left(
\begin{array}{c}
x_i' \\
x_i' + s\left(\frac{x_i'}{\alpha}\right)'
\end{array}
\right)
}
\vspace{0.2truecm}
\\
& = & 
{\displaystyle 
\frac{W_{i2}}{W_{12}} \, df_1 + \frac{W_{1i}}{W_{12}} \, df_2. 
}
\end{array}
$$
\QED

\ber
{\rm 
Note that 
$\dfrac{W_{i2}}{W_{12}}$ and $\dfrac{W_{1i}}{W_{12}}$ are $C^\infty$ functions on $t$ of order 
$a_i - a_1, a_i - a_2$ respectively. 
The above formula gives the Grassmannian lifting $\widetilde{f} : (\R^2, 0) \to \Gr(2, T\R^{N+1})$ 
of $f = \Tan(\gamma)$. 
}
\enr

\ber
\label{opening of diff}
{\rm 
If we set $g : (\R^2, 0) \to (\R^2, 0), \ g(s, t) = (f_1(s, t), f_2(s, t))$. 
Then we have that $f_3, \dots, f_{n+1} \in {\mathcal R}_g$ and that 
$f$ is an opening of $g$ in the sense of \S \ref{Opening of differentiable map-germs.}. 
}
\enr


\section{Codimension formulae and the genericity.}
\label{Codimension formulae and transversality theorems.}

We consider the jet space $J^r(I, \R P^{N+1})$. 
Let $\AAA = (a_1, a_2, \dots, a_N, a_{N+1})$ 
be a strictly increasing sequence of positive integers. 
For $r > a_{N+1}$, we define 
$$
\Sigma(\AAA) = \{ j^r\gamma(t_0) \mid t_0 \in I, \gamma : (I, t_0) \to \R P^{N+1} 
{\mbox{\rm \ is of type\ }} \AAA{\mbox{\rm  \ at\ }} t_0\}. 
$$

\bet
{\rm (\cite{Scherbak})}
$\Sigma(\AAA)$ is a semi-algebraic submanifold of codimension $\sum_{i=1}^{N+1}(a_i - i)$ in the jet space
$J^r(I, \R P^{N+1})$. 
\ent

\Proof
Let $J^r(1, N+1)$ be the fibre of the projection $\pi : J^r(I, \R P^{N+1}) \to I\times \R P^{N+1}$. 
Then $J^r(1, N+1)$ is identified with the space $\R^{(N+1)r}$ of $(N+1)\times r$-matrices. 
Then there exists an affine subspace $\Lambda \subset \R^{(N+1)r}$ 
such that $\Sigma(\AAA)$ is an image of the polynomial embedding 
$\GL(N+1, \R)\times \Lambda \to \R^{(N+1)r}$ defined by $(A, W) \mapsto AW$ for 
$A \in \GL(N+1, \R), W \in \Lambda$. Therefore $\Sigma(\AAA)$ is a semi-algebraic manifold. 

The codimension of the set consisting of jets with $\rank(W_{a_1 -1}) = 0$ is equal to 
$(N+1)(a_1 - 1)$. The codimension of the set consisting of jets with $\rank(W_{a_1 -1}) = 0, \rank(W_{a_1}) = 1$ 
and $\rank(W_{a_2-1}) = 1$ is equal to 
$(N+1)(a_1 - 1) + N(a_2 - a_1 - 1)$. Thus we have that the codimension of $\Sigma(\AAA)$ is calculated as 
$$
(N+1)(a_1 - 1)+N(a_2 - a_1 - 1) + (N-1)(a_3 - a_2 - 1) + \cdots + (a_{N+1} - a_N - 1), 
$$
which is equal to $\sum_{i=1}^{N+1}(a_i - i)$. 
\QED

\bec
For a generic curve $\gamma : I \to \R P^{N+1}$, and for any $t_0 \in I$, 
the type of $\gamma$ at $t_0$ is equal to 
$$
(1, 2, 3, \dots, N, N+1) \ {\mbox{\rm \ or \ }} (1, 2, 3, \dots, N, N+2). 
$$
\enc

\Proof
By the transversality theorem, there exists an open dense subset 
${\mathcal O} \subset C^\infty(I, \R P^{N+1})$ in $C^\infty$-topology such that 
for any $\gamma \in {\mathcal O}$ and for any $t_0 \in I$, the type $\AAA$ of 
$\gamma$ at $t_0$ satisfies $\sum_{i=1}^{N+1}(a_i - i) \leq 1$. 
Then we have $a_i = i, 1 \leq i \leq N$ and $a_{N+1} = N+1$ or $a_{N+1} = N+2$, 
and thus we have the required result. 
\QED

To treat osculating-framed curves, 
we consider the jet space of ${\mathcal C}$-integral curves, ${\mathcal C} = {\mathcal C}_{1,2,\dots,N+1}$, 
$J^r_{\mathcal C}(I,  {\mathcal F}) \subset J^r(I, {\mathcal F})$.  
Define 
$$
\Sigma_{\mathcal C}({\AAA}) := 
\{ j^r\Gamma(t_0) \mid \Gamma : (\R, t_0) \to {\mathcal F} {\mbox{\rm \ is }} {\mathcal C}{\mbox{\rm -integral, \ }} 
\pi_1\circ \Gamma {\mbox{\rm \ is of type }} 
\AAA \} 
$$ 
in $J^r_{\mathcal C}(I,  {\mathcal F})$ for sufficiently large $r$. 

\bet
{\rm (\cite{Ishikawa12}) }
\label{codim-C}
$J^r_{\mathcal C}(I,  {\mathcal F})$ is a submanifold of $J^r(I,  {\mathcal F})$ and 
the codimension of $\Sigma_{\mathcal C}({\AAA})$ in $J^r_{\mathcal C}(I,  {\mathcal F})$ is equal to 
$a_{N+1} - (N+1)$. 
\ent
\ber
{\rm 
Since any curve of finite type lifts to an ${\mathcal C}$-integral curve, 
$\Sigma_{\mathcal C}({\AAA})$ is not empty for any $\AAA$. 
}
\enr

By the transversality theorem for ${\mathcal C}$-integral curves, we have the following result: 

\bet
\label{C-types}
For a generic ${\mathcal C}$-integral curve $\Gamma : I \to {\mathcal F}_{1,2,\dots,N+1}(V^{N+2})$, 
the type $\AAA$ of $\pi_1\circ\Gamma$ 
at any point of $I$ is given by 
one of the following:
$$
\AAA = 
(1, 2, 3, \dots, N, N+1), \  \  (1, 2, \dots, i, i+2, \dots, N+1, N+2), (i = 0, \dots, N). 
$$
\ent
\Proof
By Theorem \ref{codim-C}, 
for a genetic $\Gamma$, the type of $\pi_1\circ \Gamma$ at a point in $I$ 
satisfies that $a_{N+1} - (N+1) \leq 1$, namely that 
$a_{N+1} \leq N+2$. Then we have the list of types. 
\QED

\

In general, we consider the canonical system ${\mathcal C} = {\mathcal C}_{1,2,\dots,k+1}$ 
on ${\mathcal F} = {\mathcal F}_{1,2,\dots,k+1}(V^{N+2})$, 
we consider the jet space of ${\mathcal C}$-integral curves, 
$J^r_{\mathcal C}(I,  {\mathcal F}) \subset J^r(I, {\mathcal F})$.  
Define 
$$
\Sigma_{\mathcal C}({\AAA}) := 
\{ j^rc(t_0) \mid c : (\R, t_0) \to {\mathcal F} {\mbox{\rm \ is }} {\mathcal C}{\mbox{\rm -integral, \ }} 
\pi_1\circ c {\mbox{\rm \ is of type }} 
\AAA \} 
$$ 
in $J^r_{\mathcal C}(I,  {\mathcal F})$ for sufficiently large $r$. 

\bet
\label{codim-V}
$J^r_{\mathcal C}(I,  {\mathcal F})$ is a submanifold of $J^r(I,  {\mathcal F})$ and 
the codimension of $\Sigma_{\mathcal C}({\AAA})$ in $J^r_{\mathcal C}(I,  {\mathcal F})$ 
is equal to 
$$
\sum_{i=k}^{N+1} (a_i - i) - (N - k + 1)(a_k - k). 
$$
\ent

Note that, if $k = N$, the formula is reduced to $a_{N+1} - (N+1)$ (Theorem \ref{codim-C}). 

\

\noindent
{\it Proof of Theorem \ref{codim-V}:}
Recall that ${\mathcal C} = {\mathcal C}_{1,2,\dots,k+1}$ is defined by 
$$
dx_j^{\ i-1} - x_j^{\ i} dx_{i}^{\ i-1} = 0, \quad (1 \leq i \leq k, i+1 \leq j \leq N+1)
$$
for the system of local coordinates $\left( x_j^{\ i}\right)_{0 \leq i \leq k, i+1 \leq j \leq N+1}$ of 
${\mathcal F}_{1,2,\dots,k+1}(V^{N+2})$ (\S \ref{Differential systems on flag manifolds.}). 
Then a ${\mathcal C}$-integral curve $\Gamma : I \to {\mathcal F}$ is obtained just form 
$x_{i}^{\ i-1}$-components, $1 \leq i \leq k$, and $x_j^{\ k}$-components, by integration. 
Then we see, at each point $t_0 \in I$,  $\ord(x_j^{\ 0} = \sum_{\ell = 1}^j \ord(x_j^{\ j-1})$. 
We have that the type of $\Gamma$ at $t_0$ is equal to $\AAA = (a_1, \dots, a_{N+1})$ 
if and only if 
$$
\ord(x_1^{\ 0}) = a_1, \ \ord(x_2^{\ 1}) = a_2 - a_1, \ \dots, \ \ord(x_k^{\ k-1}) = a_k - a_{k-1}, 
$$
and the type of the curve $(x_{k+1}^{\ k}, \dots, x_{N+1}^{\ k}) : (I, t_0) \to 
\R^{N-k}$ is of type $(a_{k+1} - a_k, \dots, a_{N+1} - a_k)$. 
Thus the codimension of $\Sigma_{\mathcal C}({\AAA})$ is calculated as 
$$
(a_1 - 1) + (a_2 - a_1 - 1) + \cdots + (a_k - a_{k-1} - 1) + 
\sum_{k+1}^{N+1}\left( a_j - a_k - (j - k) \right) 
= \sum_{i=k}^{N+1} (a_i - i) - (N - k + 1)(a_k - k). 
$$
\QED

\ber
Let $\pi : {\mathcal F}_{1,2,\dots,k,k+1} \to {\mathcal F}_{1,2,\dots,k}$ be the canonical projection 
defined by 
$$
\pi(V_1, V_2, \dots, V_k, V_{k+1}) = (V_1, V_2, \dots, V_k). 
$$
Then the $\pi$-fibres are projective subspaces of the flag manifold ${\mathcal F}_{1,2,\dots,k+1}$. 
In the above proof, the functions $x_{k+1}^{\ k}, \dots, x_{N+1}^{\ k}$ form a system of local projective 
coordinates of the $\pi$-fibre. 
\enr

By the transversality theorem for ${\mathcal C}$-integral curves, we have the following results:

\bet
\label{tangent-frame-list}
For a generic ${\mathcal C}_{1,2}$-integral curve $c$, the type $\AAA$ of the tangent-framed curve 
$\pi_1\circ c$ 
at any point of $I$ is given by 
one of the following: 
$$
(1, 2, 3, \dots, N, N+1), \  (1, 2, 3, \dots, N, N+2), \  (2, 3, 4, \dots, N+1, N+2). 
$$
\ent
\Proof
By Theorem \ref{codim-V}, 
for a genetic $c$, the type of $\pi_1\circ c$ at a point in $I$ 
satisfies that 
$\sum_{i=1}^{N+1} (a_i - i) - N(a_1 - 1) \leq 1$, namely that 
$\sum_{i=1}^{N+1} (a_i - i) \leq N(a_1 - 1) + 1$. Then 
$(N+1)(a_1 - 1) \leq \sum_{i=1}^{N+1} (a_i - i) \leq N(a_1 - 1) + 1$. Therefore 
$a_1 \leq 2$ and, if $a_1 = 2$, then $\AAA = (2, 3, 4, \dots, N+1, N+2)$. 
If $a_1 = 1$, then $\sum_{i=1}^{N+1} (a_i - i) \leq 1$. 
Therefore we have the result. 
\QED

\bet
\label{normal-frame-list}
For a generic ${\mathcal C}_{1,2,3}$-integral curve $\kappa$, the type $\AAA$ of 
the tangent-principal-normal-framed curve $\pi_1\circ \kappa$ 
at any point of $I$ is given by 
one of the following: 
$$
(1, 2, 3, \dots, N, N+1),  (1, 2, 3, \dots, N, N+2),  (1, 3, 4, \dots, N+1, N+2),   (2, 3, 4, \dots, N+1, N+2). 
$$
\ent

\Proof
By Theorem \ref{codim-V}, 
for a genetic $c$, the type of $\pi_1\circ c$ at a point in $I$ 
satisfies that 
$\sum_{i=2}^{N+1} (a_i - i) - (N-1)(a_2 - 2) \leq 1$, namely that 
$\sum_{i=2}^{N+1} (a_i - i) \leq (N-1)(a_2 - 2) + 1$. 
Then $N(a_2 - 2) \leq \sum_{i=2}^{N+1} (a_i - i) \leq (N-1)(a_2 - 2) + 1$, 
and we have $a_2 \leq 3$. 
If $a_2 = 3$, then $\AAA = (1, 3, 4, \dots, N+1, N+2)$ or $(2, 3, 4, \dots, N+1, N+2)$. 
If $a_2 = 2$, then $\AAA = (1, 2, 3, \dots, N, N+1)$ or $(1, 2, 3, \dots, N, N+2)$. 
\QED

\ber
{\rm 
We observe that, in all lists of 
the generic classifications of types, 
there are just three possibilities of the leading two digits: $(1, 2), (1, 3)$ and $(2, 3)$. 
These cases correspond to the cases where the projection of the tangent variety to the osculating plane 
is diffeomorphic to the map-germ $(\R^2, 0) \to (\R^2, 0)$, 
the fold singularities $(x, u) \mapsto (x^2, u)$, \lq beak to beak'
$(x, u) \mapsto (x^3 + ux^2, u)$ and 
Whitney's cusp map $(x, u) \mapsto (x^3 + ux, u)$ respectively. 
}
\enr

\section{Opening procedure of differentiable map-germs.}
\label{Opening of differentiable map-germs.}

To describe singularities of frontal mappings, we introduce the notion of \lq\lq openings" of mappings. 

The tangent variety to a curve in $\R P^{N+1}$ projects locally to the tangent variety to a space curve 
in the osculating $3$-space, and to a plane curve in the osculating $2$-plane. 
Then the tangent variety in $\R P^{N+1}$ 
can be regarded as an \lq\lq opening" of a tangent variety to a space curve and to a plane curve. 
For example, the open swallowtail, which is an opening of the swallowtail, appears in many context. 
It appears as a singular Lagrangian variety \cite{Arnol'd1}, and as 
a singular solution to certain partial differential equation \cite{Givental'}. 
The open folded umbrella appears as a \lq frontal-symplectic singularity' (\cite{IJ}).

\begin{figure}[htbp]
\begin{center}
 \includegraphics[width=7truecm, height=2.5truecm, clip, bb=51.0391 646.054 534.445 820.292]{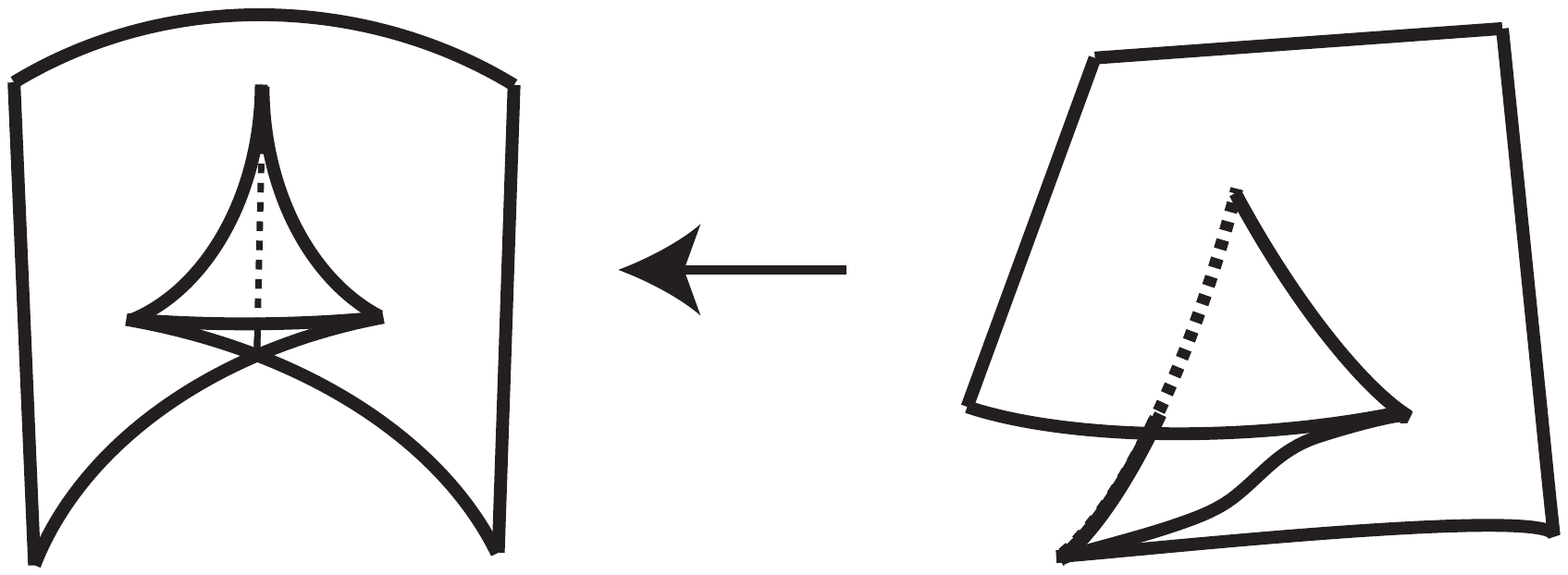}  
\end{center}
\caption{Opening of swallowtail.}
\label{Opening of swallowtail.}
\end{figure}

We denote by ${\mathcal E}_a$ the $\R$-algebra of $C^\infty$ function-germs on $(\R^n, a)$ 
with the maximal ideal ${\mathfrak m}_a$. If $a$ is the origin, then we use ${\mathcal E}_n, {\mathfrak m}_n$ 
instead of ${\mathcal E}_a, {\mathfrak m}_a$ respectively. 

\bef
{\rm 
(\cite{Ishikawa92}\cite{Ishikawa96}) 
Let $f : (\R^n, a) \to (\R^m, b)$ be a $C^\infty$ map-germ with $n \leq m$. 
We define the {\it Jacobi module} of $f$: 
$$
{\mathcal J}_f = \{ \ \sum_{j=1}^m \ p_j \ df_j \ \mid \ a_j \in {\mathcal E}_a, (1 \leq j \leq m) \ \} 
\subset \Omega_a^1, 
$$
in the space $\Omega_a^1$ of $1$-form germs on $(\R^n, a)$. 
Further we define the {\it ramification module} ${\mathcal R}_f$ by 
$$
{\mathcal R}_f := \{ h \in {\mathcal E}_a \mid dh \in {\mathcal J}_f \}. 
$$
}
\enf

Note that ${\mathcal J}_f$ is just the first order component of the graded differential ideal ${\mathcal J}^{\bullet}_f$ in $\Omega_a^{\bullet}$ 
generated by $df_1, \dots, df_m$. 
Then the singular locus is given by 
$
\Sigma_f = \{ x \in (\R^n, a) \mid \rank{\mathcal J}_f(x) < n \}. 
$
Also we consider the {\it Kernel field} 
$
\Ker(f_* : T\R^n \to T\R^m), 
$
of $f$ near $a$. 
Then we see that, for another map-germ $f' : (\R^n, a) \to (\R^{m'}, b')$ with 
${\mathcal J}_{f'} = {\mathcal J}_{f}$, $n \leq m'$, we have 
$
\Sigma_{f'} = \Sigma_{f}
$
and $\Ker(f'_*) = \Ker(f_*)$. 

Related notion was introduced in \cite{Mond3}. 

\bel
Let $f : (\R^n, a) \to (\R^m, b)$ be a $C^\infty$ map-germ. 

{\rm (1)} 
 $f^*{\mathcal E}_b \subset {\mathcal R}_f \subset {\mathcal E}_a$ and 
${\mathcal R}_f$ is an ${\mathcal E}_b$-module via $f^*$. 

{\rm (2)}
For another map-germ $f' : (\R^n, a) \to (\R^{m'}, b')$, 
${\mathcal J}_{f'} = {\mathcal J}_f$ if and only if ${\mathcal R}_{f'} = {\mathcal R}_f$. 

{\rm (3)}
If $\tau : (\R^m, b) \to (\R^m, b')$ is a diffeomorphism-germ, then 
${\mathcal R}_{\tau\circ f} = {\mathcal R}_f$. If $\sigma : 
(\R^n, a') \to (\R^n, a)$ is a diffeomorphism-germ, then 
${\mathcal R}_{f\circ\sigma} = \sigma^*({\mathcal R}_f)$. 
\enl

\Proof (1) follows from that,  
if $h \in {\mathcal R}_f$ and $dh = \sum_{j=1}^m p_j df_j$, then we have
$$
d\{(k\circ f)h\} = \sum_{j=1}^m\left\{ (k\circ f)p_j + h\left({\pa k}/{\pa y_j}\right)\right\} df_j. 
$$

(2) It is clear that ${\mathcal J}_{f'} = {\mathcal J}_f$ implies 
${\mathcal R}_{f'} = {\mathcal R}_f$. Conversely suppose ${\mathcal R}_{f'} = {\mathcal R}_f$. 
Then any component $f'_j$ of $f'$ belongs to 
${\mathcal R}_{f'} = {\mathcal R}_f$, hence $df_j \in {\mathcal J}_f$. 
Therefore ${\mathcal J}_{f'} \subset {\mathcal J}_f$. 
By the symmetry we have ${\mathcal J}_{f'} = {\mathcal J}_f$. 

(3) follows from that 
${\mathcal J}_{\tau\circ f} = {\mathcal J}_f$ and 
${\mathcal J}_{f\circ\sigma} = \sigma^*({\mathcal J}_f)$. 
\QED

\bef
{\rm 
Let $f : (\R^n, a) \to (\R^m, b), \ n \leq m$ be a $C^\infty$ map-germ. 
Given $h_1, \dots, h_r \in {\mathcal R}_f$, the map-germ 
$F : (\R^n, 0) \to \R^m\times\R^r = \R^{m+r}$ defined by 
$$
f = (f_1, \dots, f_m, h_1, \dots, h_r)
$$
is called an {\it opening} of $f$, while $f$ is called a {\it closing} of $F$. 

An opening $F = (f, h_1, \dots, h_r)$ of $f$ is called a {\it versal opening} (resp. 
{\it mini-versal opening}) of $f$, if 
$1, h_1, \dots, h_r$ form a (minimal) system of 
generators of ${\mathcal R}_f$ as an ${\mathcal E}_b$-module via $f^*$. 
}
\enf

Note that an opening of an opening of $f$ is an opening of $f$. 

\

Here we summarise known results on the ramification module. 
A map-germ $f : (\R^n, a) \to (\R^m, b)$ is called {\it finite} if 
$\dim_{\R}{\mathcal E}_a/(f^*{\mathfrak m}_b){\mathcal E}_a < \infty$. 

\bep
\label{finiteness}
{\rm (Theorem 3.1 of \cite{Ishikawa94}, Corollary 2.4 of \cite{Ishikawa96})}
If $f : (\R^n, a) \to (\R^m, b)$ 
is finite and of corank at most one. 
Then we have 

{\rm (1)} 
${\mathcal R}_f$ is a finite ${\mathcal E}_b$-module. 
Therefore there exists a versal opening of $f$. 

{\rm (2)}
$1, h_1, \dots, h_r \in {\mathcal R}_f$ generate ${\mathcal R}_f$ 
as ${\mathcal E}_b$-module if and only if 
they generate the vector space ${\mathcal R}_f/f^*({\mathfrak m}_b){\mathcal R}_f$ over $\R$. 
\enp

\ber
\label{basis}
{\rm 
By Proposition \ref{finiteness}, we see that 
$1, h_1, \dots, h_r \in {\mathcal R}_f$ form a minimal system of generators of ${\mathcal R}_f$ 
as ${\mathcal E}_b$-module if and only if 
they form a basis of $\R$-vector space ${\mathcal R}_f/f^*({\mathfrak m}_b){\mathcal R}_f$. 
}
\enr

Let $k \geq 0, m \geq 0$.  To present the normal forms of Morn map, 
consider variables 
$t, \lambda = (\lambda_1, \dots, \lambda_{k-1}), \mu = \left( \mu_{ij}\right)_{1\leq i \leq m, 1 \leq  j \leq k}$ 
and polynomials 
$$
F(t, \lambda) = t^{k+1} + \sum_{i=1}^{k-1} \lambda_jt^j, \quad 
G_i(t, \mu) = \sum_{j=1}^k \mu_{ij}t^j, (1 \leq i \leq m). 
$$
Let 
$f : (\R^{k + km}, 0) 
\to (\R^{m + k + km}, 0)$ be a Morin map
defined by 
$$
f(t, \lambda, \mu) := (F(t, \lambda), G(t, \mu), \lambda, \mu), 
$$
for the above polynomials $F$ and $G$. 

For $\ell \geq 0$, we denote by $F_{(\ell)}, G_{i(\ell)}$ the polynomials 
$$
F_{(\ell)}(t, \lambda) = 
\dint{0}{t} s^\ell F(s, \lambda) ds, \quad 
G_{i\, {(\ell)}}(t, \mu) = \dint{0}{t} s^\ell G_i(s, \mu) ds. 
$$
Then we have: 

\bep
\label{generator}
{\rm (Theorem 3 of \cite{Ishikawa92})} 
The ramification module ${\mathcal R}_f$ of the Morin map $f$ is minimally generated over 
$f^*{\mathcal E}_{m + k + km}$ by the $1+k+(k-1)m$ elements 
$$
1, \ F_{(1)}, \ \dots, \ F_{(k)}, \ G_{1\, (1)}, \ \dots, \ G_{1\, (k-1)}, \ \dots, \ 
G_{m\, (1)}, \ \dots, \ G_{m\, (k-1)}. 
$$
The map-germ ${\mathbf F} : (\R^{k+mk}, 0) \to (\R^{m+k+km}\times \R^{k+(k-1)m}, 0) 
= (\R^{2(k+km)}, 0)$ defined by 
$$
{\mathbf F} = \left( f, \ F_{(1)}, \dots, F_{(k)}, G_{1\, (1)}, \dots, G_{1\, (k-1)}, \dots, 
G_{m\, (1)}, \dots, G_{m\, (k-1)}\right) 
$$
is a versal opening of $f$. 
\enp

\ber
{\rm 
It is shown in \cite{Ishikawa92} moreover 
that ${\mathbf F}$ is an isotropic map for a symplectic structure 
on $\R^{2(k+km)}$. 
}
\enr

\bep
\label{generators of ramification}
{\rm (cf. Proposition 1.6 of \cite{Ishikawa94}, Lemma 2.4 of \cite{Ishikawa95})} 
Let $f : (\R^n, a) \to (\R^m, b)$ be a $C^\infty$ map-germ and 
$
F : (\R^{n+\ell}, (a,0)) \to (\R^{m+\ell}, (b,0))
$ be an unfolding of $f$: 
$F(x, u) = (F_1(x, u), u)$ and $F_1(x, 0) = f(x)$. 
Let $i : (\R^n, a) \to (\R^{n+\ell}, (a,0))$ be the inclusion, $i(x) = (x, 0)$. 
Then we have:

{\rm (1)}  
$i^*{\mathcal R}_F \subset {\mathcal R}_f$. 

{\rm (2)} 
If $f$ is of corank $\leq 1$ with $n \leq m$, then 
$i^*{\mathcal R}_F = {\mathcal R}_f$. 
If $1, H_1, \dots, H_r$ generate ${\mathcal R}_F$ via $F^*$, then 
$1, i^*H_1, \dots, i^*H_r$ generate ${\mathcal R}_f$ via $f^*$. 

{\rm (3)} Let $\ell$ be a positive integer and  $F = (F_1(t, u), u) : (\R^n, 0) \to (\R^n, 0)$ an unfolding 
of $f : (\R, 0) \to (\R, 0), f(t) = F_1(t, 0) = t^{\ell}$. 
Suppose $H_1, \dots, H_r \in {\mathcal R}_F \cap {\mathfrak m}_n$. 
Then $1, H_1, \dots, H_r$ generate ${\mathcal R}_F$ via $F^*$ if and only 
$i^*H_1, \dots, i^*H_r$ generate ${\mathfrak m}_1^{\ell+1}/{\mathfrak m}_1^{2\ell}$. 
In particular 
$F_{1(1)}, \dots, F_{1(\ell-1)}$ form a system of generators of ${\mathcal R}_F$ via $F^*$ over 
${\mathcal E}_n$. 
\enp

\Proof
(1) is clear. (2) 
Let $H \in {\mathcal R}_F$. Then $dH \in {\mathcal J}_F$. Hence 
$d(i^*H) = i^*(dH) \in i^*{\mathcal J}_F \subset {\mathcal J}_f$. Therefore $i^*H \in {\mathcal R}_f$. 
Let $f$ be of corank at most one. Suppose $h \in {\mathcal R}_f$. Then 
$dh = \sum_{j=1}^m a_j df_j$ for some $a_j \in {\mathcal E}_a$. There exist
$A_j, B_k \in {\mathcal E}_{(a,0)}$ such that $i^*A_j = a_j$ and 
the $1$-form $\sum_{j=1}^m A_j d(F_1)_j + \sum_{k=1}^\ell B_k d\lambda_k$ is closed (cf. Lemma 2.5 of 
\cite{Ishikawa96}). Then there exists an $H \in {\mathcal E}_{(a,0)}$ such that 
$dH = \sum_{j=1}^m A_j d(F_1)_j + \sum_{k=1}^\ell B_k d\lambda_k \in {\mathcal J}_F$ and 
$d(i^*H) = i^*(dH) = dh$. Then there exists $c \in \R$ such that $h = i^*H + c = i^*(H + c)$, 
and $H + c \in {\mathcal R}_F$.  Therefore $h \in i^*{\mathcal R}_F$. 
Since $i^*$ is a homomorphism over $j^* : {\mathcal E}_{(b,0)} \to {\mathcal E}_b$, where 
$j : (\R^m, 0) \to (\R^{m+\ell}, 0)$ is the inclusion $j(y) = (y, 0)$, we have the consequence. 
(3) It is easy to show that ${\mathcal R}_f = \R + {\mathfrak m}_1^{\ell}$. 
By Proposition \ref{finiteness} (2), $1, H_1, \dots, H_r$ generate ${\mathcal R}_F$ as 
${\mathcal E}_n$-module via $F^*$ if and only if 
they generate ${\mathcal R}_F/F^*({\mathfrak m}_n){\mathcal R}_F$ over $\R$. 
Since 
$$
{\mathcal R}_F/F^*({\mathfrak m}_n){\mathcal R}_F 
\cong (\R + {\mathfrak m}_1^{\ell})/(f^*{\mathfrak m}_1)(\R + {\mathfrak m}_1^{\ell}) 
\cong {\mathfrak m}_1^{\ell+1}/{\mathfrak m}_1^{2\ell}
$$ 
we have the consequence. 
\QED

\bep
\label{mini-versal}
Let $f : (\R^n, a) \to (\R^m, b), \ n \leq m$ be a $C^\infty$ map-germ. 
\\
{\rm (1)} 
For any versal opening $F : (\R^n, a) \to (\R^{m+r}, F(a))$ of $f$ and for any 
opening $G : (\R^n, a) \to (\R^{m+s}, G(a))$, there exists an affine bundle map 
$\Psi : (\R^{m+r}, F(a)) \to (\R^{m+s}, G(a))$ over $(\R^m, f(a))$ 
such that $G = \Psi\circ F$. 
\\
{\rm (2)} 
For any mini-versal openings $F : (\R^n, a) \to (\R^{m+r}, F(a))$ and 
$F' : (\R^n, a) \to (\R^{m+r}, F'(a))$ of $f$, there exists 
an affine bundle isomorphism 
$\Phi : (\R^{m+r}, F(a)) \to  (\R^{m+r}, F'(a))$ over $(\R^m, f(a))$ 
such that $F' = \Psi\circ F$. 
In particular, the diffeomorphism class of mini-versal opening of $f$ is unique. 
\\
{\rm (3)} 
Any versal openings $F'' : (\R^n, a) \to (\R^{m+s}, F''(a))$ 
of $f$ is diffeomorphic to $(F, 0)$ for a mini-versal opening of $f$. 
\enp

\Proof
(1) Let $F = (f, h_1, \dots, h_r)$ and $G = (f, k_1, \dots, k_s)$. 
Since $k_j \in {\mathcal R}_f$, there exist $c_j^{\ 0}, c_j^{\ 1}, \dots, c_j^{\ r} \in {\mathcal E}_b$ 
such that $k_j = c_j^{\ 0}\circ f + (c_j^{\ 1}\circ f)h_1 + \cdots +  (c_j^{\ r}\circ f)h_r$. 
Then it suffices to set $\Psi(y, z) = (y, (c_j^{\ 0}(y) + c_j^{\ 1}(y)z_1 + \cdots +  c_j^{\ r}(y)z_r)_{1\leq j\leq s})$. 
(2) By (1) there exists an affine bundle map $\Psi$ with $F' = \Psi\circ F$. 
From the minimality, we have that the matrix $(c_j^{\ i}(b))$ is regular. (See Remark \ref{basis}). 
Therefore $\Psi$ is a diffeomorphism-germ. 
(3) Let $F = \Psi\circ F''$ for some affine bundle map $\Psi$. Then the matrix $(c_j^{\ i}(b))$ is of rank $r$. 
Therefore $F''$ is diffeomorphic to $(F, k_1, \dots, k_{s-r})$ for some 
$k_j \in {\mathcal R}_f$. Write each $k_j = K_j\circ F$ for some 
$K_j \in {\mathcal E}_{F(a)}$. Then we set $\Xi(y, z, w) = (y, z, w - K\circ F)$. Then $\Xi$ is a 
local diffeomorphism on $\R^{m+r+(s-r)}$ and $\Xi\circ(F, k_1, \dots, k_{s-r}) = (F, 0)$. 
\QED

\section{Normal forms of tangent surfaces.}
\label{Normal forms of tangent surfaces.}

According to a geometric restriction expressed in differential system, 
we have imposed on curves in projective spaces a system of differential 
equations (\S \ref{Differential systems on flag manifolds.}). 
The genericity, in such a restricted class of curves, naturally 
implies a restriction on types of curves (\S \ref{Codimension formulae and transversality theorems.}). 
Then we use the following results to solve the classification problem. 
For the concrete expression of normal forms, see \S \ref{Differential systems on flag manifolds.}. 

\bet 
\label{normal forms}
{\rm (1)} In $\R P^3$, 
the tangent variety of a curve of type $(1, 2, 3)$ {\rm (}resp. $(1, 2, 4)$, $(2, 3, 4)$, $(1, 3, 4)${\rm )} 
is locally diffeomorphic to the cuspidal edge {\rm (}the folded umbrella, the swallowtail, the Mond surface{\rm )} 
in $\R^3$.  

{\rm (2)} {\rm (Higher codimensional case.)} 
In $\R P^{N+1}, N+1 \geq 4$, 

{\rm (i)} the tangent variety of a curve of type $(1, 2, 3, a_4, \dots, a_{N+1})$ 
is locally diffeomorphic to the cuspidal edge 
$(\R^2, 0) \to (\R^3, 0)$ composed with the inclusion to $(\R^{N+1}, 0)$. 

{\rm (ii)} the tangent variety of a curve of type 
$(1, 3, 4, 5, a_5, \dots, a_{N+1})$ 
is locally diffeomorphic to 
the open Mond surface $(\R^2, 0) \to (\R^4, 0)$ composed with the inclusion to $(\R^{N+1}, 0)$. 

{\rm (iii)} the tangent variety of a curve of type 
$(2, 3, 4, 5, a_5, \dots, a_{N+1})$
is locally diffeomorphic to 
the open swallowtail $(\R^2, 0) \to (\R^4, 0)$ composed with the inclusion to $(\R^{N+1}, 0)$. 
\ent

\Proof
(1) is proved in Theorem 1 $(n = 2)$ in \cite{Ishikawa95}. 
(2) In each case, the idea is to show 
that the tangent map-germ $\Tan(\gamma)$ is diffeomorphic to 
a mini-versal opening of an appropriate map-germ: 

(i) the fold map-germ $(\R^2, 0) \to (\R^2, 0)$. 

(ii) the Mond surface $(\R^2, 0) \to (\R^3, 0)$. 

(iii) the swallowtail $(\R^2, 0) \to (\R^3, 0)$. 

Then, by Proposition \ref{mini-versal}, the diffeomorphism class of 
the tangent map-germ is unique and we get the required results. 

Let $\gamma : (\R, 0) \to \R P^{N+1}$ 
be a curve-germ of type $(a_1, a_2, \dots, a_{N+1})$,   
$$
\gamma(t) = \left( x_1(t), x_2(t), \dots, x_{N+1}(t)\right)
$$
a local affine representation of $\gamma$ as in Lemma \ref{affine coord}, and 
$$
f(s, t) = (f_1(s, t), f_2(s, t), \dots, f_{N+1}(s, t)) = 
\left( x_i(t) + s\, \frac{1}{\alpha(t)}x_i'(t)\right)_{1 \leq i \leq N+1}, 
$$
the parametrisation of the tangent variety to $\gamma$, 
where $\alpha(t) = t^{a_1 - 1}$. We may suppose $x_1(t) = t^{a_1}$. 

We define $g' : (\R^2, 0) \to (\R^2, 0)$ by $g' = (f_1, f_2)$.
Then, by Lemma \ref{Wronskian} and Remark \ref{opening of diff}, we see 
that 
$f_3, \dots, f_{N+1} \in {\mathcal R}_{g'}$. 
Note that $f_1(s, t) = x_1(t) + a_1 s$ is a regular function. We regard $f_1(s, t)$ as an unfolding 
parameter $u$. 
Then there exist diffeomorphism-germ $\sigma : (\R^2, 0) \to (\R^2, 0)$ and 
$\tau : (\R^2, 0) \to (\R^2, 0)$ such that $\sigma$ is of form $\sigma(u, t) = (\sigma_1(u), t\sigma_2(u, t))$ 
and $g = \tau\circ g'\circ\sigma$ is equal to 
(i) $(u, t) \mapsto (u, t^2 + ut)$, (ii) $(u, t) \mapsto (u, t^3 + ut^2)$, (iii) $(u, t) \mapsto (u, t^3 + ut)$. 
Then $f_3\circ\sigma, \dots, f_{N+1}\circ\sigma$ belongs ${\mathcal R}_{g} 
= {\mathcal R}_{g'\circ\sigma}$. 
Then, by Lemma \ref{generators of ramification}, 
(i) $F = (f_1\circ \sigma, f_2\circ\sigma, f_3\circ\sigma)$, (ii)(iii) $F = (f_1\circ \sigma, f_2\circ\sigma, f_3\circ\sigma,  f_4\circ\sigma)$, are versal opening of $g$ respectively. 
Note that in cases (ii) and (iii), $F$ is a versal opening of also Mond surface and swallowtail respectively.  
Then, by Proposition \ref{generators of ramification} (3), 
we have that $f\circ\sigma$ is diffeomorphic to 
(i) $(u,  t^2 + ut, \frac{2}{3}t^3 + \frac{1}{2}ut^2, 0, \dots, 0)$, (ii) $(u, t^3 + ut^2, \frac{3}{4}t^4 + \frac{2}{3}ut^3, 
\frac{3}{5}t^5 + \frac{1}{2}ut^4, 0, \dots, 0)$, (iii) $(u, t^3 + ut^2, \frac{3}{4}t^4 + \frac{1}{2}ut^2, 
\frac{3}{5}t^5 + \frac{1}{3}ut^3, 0, \dots, 0)$, as required. 
\QED

\bet
\label{1,2,4,5}
The tangent variety of a curve of type 
$(1, 2, 4, 5, a_5, \dots, a_{N+1})$ 
is locally diffeomorphic to 
the open folded umbrella $(\R^2, 0) \to (\R^4, 0)$ composed with the inclusion to $(\R^{N+1}, 0)$. 
\ent

\Proof 
We argue as in Theorem \ref{normal forms}. However in this case 
$(f_1\circ \sigma, f_2\circ\sigma, f_3\circ\sigma, f_4\circ\sigma)$ 
is not a versal opening of $g = (u, t^2 + ut)$. 
(In fact the open folded umbrella is not a versal opening of the folded umbrella $(\R^2, 0) \to (\R^3, 0)$. )

To show Theorem \ref{1,2,4,5}, we define 
$$
{\mathcal R}_{g}^{(2)} := \left\{ h \in t^2{\mathcal E}_2 \mid dh 
\in t^2{\mathcal J}_{g} \right\}. 
$$
Then $f_i\circ\sigma \in {\mathcal R}_{g}^{(2)}, (i \geq 3)$. 
We see that 
$f_3\circ\sigma, f_4\circ\sigma$ generate 
${\mathcal R}_g^{(2)}$ over $g^*{\mathcal E}_2$. 
In fact $h_1, \dots, h_r$ generate ${\mathcal R}_g^{(2)}$ as ${\mathcal E}_2$-module 
if and only if $i^*h_1, \dots, i^*h_r$ generate ${\mathfrak m}_1^4/{\mathfrak m}_1^6$ over $\R$. 
(See Lemma 2.4 of \cite{Ishikawa95}). 
Also $h_1 = \frac{1}{2}t^4 + \frac{1}{3}ut^3, h_2 = 
\frac{2}{5}t^5 + \frac{1}{4}ut^4$ generate ${\mathcal R}_g^{(2)}$. 
We write 
$f_i\circ\sigma = (a_i\circ g)h_1 + (b_i\circ g)h_2, (i \geq 3)$, for some $a_i, b_i \in {\mathcal E}_2$. 
We define 
$\Psi : (\R^{N+1}, 0) \to (\R^{N+1}, 0)$ by 
$$
\begin{array}{rcl}
\Psi(x) & = & (x_1, x_2, a_3(x_1, x_2)x_3 + b_3(x_1, x_2)x_4, a_4(x_1, x_2)x_3 + b_4(x_1, x_2)x_4, 
\\
& & 
\quad\quad 
x_i - a_i(x_1, x_2)x_3 + b_i(x_1, x_2)x_4 (5 \geq i)). 
\end{array}
$$
Then $\Psi$ is a diffeomorphism-germ and $\Psi\circ f\circ\sigma = (g, h_1, h_2, 0)$. 
Thus we have that 
$f\circ \sigma$ is diffeomorphic to 
$
(g, h_1, h_2, 0) = (u, t^2 + ut, \frac{1}{2}t^4 + \frac{1}{3}ut^3, \frac{2}{5}t^5 + \frac{1}{4}ut^4, 0, \dots, 0)
$ 
as required. 
\QED

\

\

\noindent
{\it Proofs of the classification theorems.}
Theorems \ref{tangent to generic curve}, \ref{tangent-framed classification}, 
\ref{tangent-principal-normal-framed classification}, \ref{osculating-frame-geq4} follows from 
Theorems \ref{tangent-frame-list}, \ref{normal-frame-list}, \ref{C-types} and 
Theorems \ref{normal forms}, \ref{1,2,4,5}.

\

We are led, in our generic classifications in a geometric setting, to find the following result, 
which we use in \S \ref{Singularities on tangent varieties to osculating framed contact-integral curves.}.

\bet
\label{1,3,4,6}
The tangent variety of a curve of type 
$(1, 3, 4, 6, a_5, \dots, a_{N+1})$ in $\R P^{N+1}, N+1 \geq 4$, 
has unique diffeomorphism class. 
\ent

We may call it the \lq unfurled Mond surface', distinguished with the open Mond surface. 
The normal form $(\R^2, 0) \to (\R^{N+1}, 0)$ of the {\it unfurled Mond surface} is given by 
$$
(s, t) \mapsto \left(t + s, \ t^3 + 3st^2, \ t^4 + 4st^3, \ t^6 + 6st^5, \ 0, \dots, 0\right), 
$$
which is diffeomorphic to 
$$
(x, u) \mapsto \left(u, \ t^3 + ut^2, \ \frac{3}{4}t^4 + \frac{2}{3}ut^3, \ \frac{1}{2}t^6 + \frac{2}{5}ut^5, \ 
0, \dots, 0\right). 
$$

To show Theorem \ref{1,3,4,6}, we prepare the following:

\bel
\label{Lemma(1,3,4,6)}
{\rm (cf. Lemma 2.4 of \cite{Ishikawa95})} 
Let $g : (\R^2, 0) \to (\R^2, 0)$ be the map-germ defined by $g(t, u) = (u, t^3 + ut^2)$. We set 
$$
{\mathcal R}_g^{(3)} := \left\{ h \in t^3{\mathcal E}_2 \mid dh \in t^3{\mathcal J}_g \right\}. 
$$
and set $T = t^3 + ut^2, T_i = \frac{3}{i+3}t^{i+3} + \frac{2}{i+2}ut^{i+2}, (i = 1, 2, 3, \dots)$. 
Then we have 
{\rm (1)}
${\mathcal R}_g^{(3)} = {\mathcal R}_g \cap t^5{\mathcal E}_2$. 
{\rm (2)}
${\mathcal R}_g^{(3)}$ is a finite ${\mathcal E}_2$-module via $g^* : {\mathcal E}_2 \to {\mathcal E}_2$ 
generated by $T_3, TT_1, T_1^2$. 
{\rm (3)} 
Let $\iota : (\R, 0) \to (\R^2, 0), \iota(t) = (t, 0)$. 
Then 
$h_1, \dots, h_\ell \in {\mathcal R}_g^{(3)}$ generate ${\mathcal R}_g^{(3)}$ as 
${\mathcal E}_2$-module via $g^*$ if and only if 
$\iota^*h_1, \dots, \iota^*h_\ell$ generate $t^6{\mathcal E}_1/t^9{\mathcal E}_1$ over $\R$. 
{\rm(}Note that $T_1 \not\in {\mathcal R}_g^{(3)}$.{\rm)}
\enl

\Proof
(1) First note that ${\mathcal R}_g^{(3)} = \left\{ h \in t^3{\mathcal E}_2 \ \left\vert\ \dfrac{\pa h}{\pa t} 
\in t^3\dfrac{\pa T}{\pa t}{\mathcal E}_2 \right\}\right.$. 
Let $h \in {\mathcal R}_g^{(3)}$. Then $\dfrac{\pa h}{\pa t} \in t^4{\mathcal E}_2$ and 
$h \in t^3{\mathcal E}_2$. Therefore $h \in {\mathcal R}_g \cap t^5{\mathcal E}_2$. Conversely 
let $h \in {\mathcal R}_g \cap t^5{\mathcal E}_2$. Then $\dfrac{\pa h}{\pa t} = t^3\dfrac{\pa T}{\pa t}K$ 
for some $K \in {\mathcal E}_2$. Since $h(0, 0) = 0$, 
we have $\dfrac{\pa h}{\pa u} \in t^5{\mathcal E}_2$. Therefore 
$dh \in t^3{\mathcal J}_g$ and $h \in {\mathcal R}_g^{(3)}$. Thus we have the equality. 

(2) Let $h \in {\mathcal R}_g^{(3)}$. 
Then $h = a\circ g + b\circ g T_1 + c\circ g T_2$, for some $a, b, c \in {\mathcal E}_2$. 
Since $h \in t^5{\mathcal E}_2$, 
$h = \widetilde{a}\circ g T^3+ 
\widetilde{b}\circ g TT_1 + \widetilde{c}\circ g TT_2$, for some 
$\widetilde{a}, \widetilde{b}, \widetilde{c} \in {\mathcal E}_2$. 
Note that $T^3, TT_1, TT_2 \in {\mathcal R}_g^{(3)}$. 
Moreover we have directly 
$$
T^3 = \dfrac{32}{15}uT_1^2 + 2TT_3 + \dfrac{14}{3}T_4, \ 
TT_2 = \dfrac{16}{15}T_1^2 + \dfrac{7}{3}uT_4,  \ 
T_4 = \dfrac{4}{7}TT_1 - \dfrac{20}{21}uT_3. 
$$
Therefore we have 
$$
TT_2 =  - \dfrac{20}{9}u^2T_3 + \dfrac{4}{3}uTT_1  + \dfrac{16}{15}T_1^2,  \ 
T^3 = \left(2T - \dfrac{40}{9}u^3\right)T_3 + \dfrac{8}{3}u^2TT_1 + \dfrac{32}{15}uT_1^2.  
$$

(3) $\iota^* : {\mathcal E}_2 \to {\mathcal E}_1$ induces 
$\iota^* : {\mathcal R}_g^{(3)} \to t^6{\mathcal E}_1$, which is clearly surjective. 
Moreover we have $(\iota^*)^{-1}(t^9{\mathcal E}_1) = g^*{\mathfrak m}_2{\mathcal R}_g^{(3)}$. 
Therefore $\iota^*$ induces an isomorphism 
${\mathcal R}_g^{(3)} / g^*{\mathfrak m}_2{\mathcal R}_g^{(3)} \cong t^6{\mathcal E}_1/t^9{\mathcal E}_1$ 
as $\R$-vector spaces. 
By (2) and by Malgrange-Mather's preparation theorem \cite{Brocker}, we have the required result. 
\QED

\

\noindent
{\it Proof of Theorem \ref{1,3,4,6}:}
We give the proof for the case $N+1 = 4$. In general case we can argue similarly. 

Let $\gamma : (\R, 0) \to \R P^4$ be a curve of type $(1, 3, 4, 6)$. 
The tangent map-germ $\Tan(\gamma)$ is an opening of a Mond surface. 
However it is not versal. So we need a specialised idea to show the determinacy result in this situation. 
Let 
$$
\gamma(t) = (t, t^3 + \varphi(t), t^4 + \psi(t), t^6 + \rho(t)), 
$$
with $\varphi \in \mathfrak{m}_1^4, \psi \in \mathfrak{m}_1^5, \rho \in \mathfrak{m}_1^7$. 
Then $f = \Tan(\gamma)$ is given by 
$$
f(s, t) = 
\left(t + s, \ t^3 + 3st^2 + \Phi(t), \ t^4 + 4st^3 + \Psi(t), \ t^6 + 6st^5 + R(t)\right), 
$$
where $\Phi(s, t) = \varphi(t) + s\varphi'(t), \Psi(s, t) = \psi(t) + s\psi'(t), R(s, t) =  \rho(t) + s\rho'(t)$. 
We set $u = t+s$. Then 
$$
f(u, t) = 
\left(u, \ -2t^3 + 3ut^2 + \widetilde{\Phi}(t), \ -3t^4 + 4ut^3 + \widetilde{\Psi}(t), \ -5t^6 + 6ut^5 + 
\widetilde{R}(t)\right), 
$$
where $\Phi(s, t) = \varphi(t) + (u - t)\varphi'(t), \Psi(s, t) = \psi(t) + (u - t)\psi'(t), R(s, t) =  \rho(t) + (u - t)\rho'(t)$. 
From the determinacy of tangent varieties to curves of type $(1, 3, 4)$ in $\R^3$ 
(\cite{Mond2}, \cite{Ishikawa93}), 
we have that there exist diffeomorphism-germ $\sigma : (\R^2, 0) \to (\R^2, 0)$ of form 
$\sigma(u, t) = (\sigma_1(u), t\sigma_2(u, t))$ and a diffeomorphism-germ  
$\tau : (\R^{4}, 0) \to (\R^{4}, 0)$ 
such that 
$$
f\circ\sigma(u, t) = \left(u, \ T(u, t), \ T_1(u, t), \ T_3(u, t) + S_3(u, t)\right), 
$$
with 
$$
T = t^3 + ut^2, \quad T_1 = \dfrac{3}{4}t^4 + \dfrac{2}{3}ut^3, \quad T_3 = \dfrac{1}{2}t^6 + \dfrac{2}{5}ut^5, 
$$
$S_3 \in {\mathcal R}_g^{(3)}$, $g = (u, \ t^3 + ut^2)$, 
$\iota^*S_3 \in {\mathfrak m}_1^7$. 
Then we have, by Lemma \ref{Lemma(1,3,4,6)}, 
$$
S_3 = A_3\circ g \ T_3 + B_3\circ g \ TT_1 + C_3\circ g \ T_1^2, 
$$
for some 
$A_3, B_3, C_3 \in {\mathcal E}_2$ with 
$A_3(0, 0) = 0$. 
Define $\Xi : (\R^4, 0) \to (\R^4, 0)$ by 
$$
\begin{array}{lcr}
\Xi(x_1, x_2, x_3, x_4) 
& = & 
\left(x_1, \ x_2, \ x_3 + A_1(x_1, x_2)x_4 + B_1(x_1, x_2)x_2x_3 + C_1(x_1, x_2)x_3^2, \right.
\vspace{0.1truecm}
\\
& & 
\left.
x_4 + A_3(x_1, x_2)x_4 + B_3(x_1, x_2)x_2x_3 + C_3(x_1, x_2)x_3^2\right). 
\end{array}
$$
Then the Jacobi matrix of $\Xi$ is the unit matrix, so $\Xi$ is a diffeomorphism-germ and 
$$
\Xi^{-1}\circ f\circ\sigma = \left(u, \ t^3 + ut^2, \ \dfrac{3}{4}t^4 + \dfrac{2}{3}ut^3, \ 
\dfrac{1}{2}t^6 + \dfrac{2}{5}ut^5\right). 
$$
\QED

\section{Singularities on tangent varieties to osculating-framed contact-integral curves.}
\label{Singularities on tangent varieties to osculating framed contact-integral curves.}

We give results on the classification of 
singularities of tangent varieties to contact-integral curves 
(resp.  osculating framed contact-integral curves) in a contact projective space. 

Let $V$ be a symplectic vector space of dimension $2n+2$.  
Consider the isotropic flag manifold: 
$$
{\mathcal F}_{\Lag} = {\mathcal F}_{\Lag}(V) := 
\{ V_1 \subset V_2 \subset \cdots \subset V_n \subset V_{n+1} \subset V 
\mid V_{n+1}  \ {\mbox {\rm is Lagrangian}} \}. 
$$
Note that ${\mathcal F}_{\Lag}$ is a finite quotient of 
${\mbox{\rm U}}(n+1), \dim ({\mathcal F}_{\Lag}) = (n+1)^2$ and that
${\mathcal F}_{\Lag}(V)$ is embedded into ${\mathcal F}(V) = {\mathcal F}_{1,2,\dots,n+1,\dots, 2n+1}(V)$ by taking symplectic orthogonals:  
$$
(V_1, V_2, \dots, V_n, V_{n+1}) \mapsto (V_1, V_2, \dots, V_n, V_{n+1}, V_n^s, \dots, V_2^s, V_1^s), 
$$

Define a differential system ${\mathcal E} \subset T{\mathcal F}_{\Lag}$ by 
$$
v \in {\mathcal E}_{(V_1, \dots, V_{n+1})} \Longleftrightarrow  
{\pi_i}_*(v) \in T\,{\Gr}(i, V_{i+1}) (\subset T\,{\IGr}(i, V)), (1 \leq i \leq n). 
$$
where $\IGr$ means the isotropic Grassmannian, 
$\pi_i : {\mathcal F}_{\Lag} \to {\IGr}(i, V)$ is the canonical projection. 
Then $\rank({\mathcal E}) = n+1$ and ${\mathcal E}$ is bracket generating. 

If $n = 1$, then we have $\dim {\mathcal F}_{\Lag} = 4$ and ${ {\mathcal E}}$ 
is an Engel structure on ${\mathcal F}_{\Lag}$ (\cite{IMT}). 

An ${\mathcal E}$-integral curve $c : I \to {\mathcal F}_{\Lag}$ is a $C^\infty$ family 
$$
(V_1(t), V_2(t), \dots, V_n(t), V_{n+1}(t))
$$ of
isotropic flags in the symplectic vector space $V$ 
such that $V_i(t)$ moves momentarily in $V_{i+1}(t)$. 

\ber
{\rm 
The projective space $P(V^{2n+2}) \cong \R P^{2n+1}$ 
has the canonical contact structure ${\mathcal D} \subset T(P(V))$ : 
For $V_1 \in P(V)$ and for $v \in T_{V_1}P(V)$, we define 
$$
v \in  {\mathcal D}_{V_1} \Longleftrightarrow \pi_{1*}(v) \in T(P(V_1^s)) (\subset T(P(V))). 
$$
If $c : I \to {\mathcal F}_{\Lag}(V)$ is an ${\mathcal E}$-integral curve, then  
$\gamma = \pi_1\circ c : I \to P(V)$ is a ${\mathcal D}$-integral curve. 
}
\enr

We consider 
the space $J^r_{\mathcal E}(I, {\mathcal F}_{\Lag}(\R^{2n+2})$ of ${\mathcal E}$-integral jets 
in 
\\
$J^r(I, {\mathcal F}_{\Lag}(\R^{2n+2}))$ and set 
$$
\Sigma_{\mathcal E}(\AAA)  := 
\{ j^r\Gamma(t_0) \mid 
t_0 \in I, \Gamma : (\R, t_0) \to {\mathcal F}_{\Lag}(\R^{2n+2}) {\mbox{\rm \ is }} 
{\mathcal E}{\mbox{\rm -integral, \ }} 
 \pi_1\circ \Gamma {\mbox{\rm \ is of type }} \AAA \}. 
$$
Then we have the codimension formula 
for osculating framed contact-integral curves. 

\bet
\label{codim-E}
The set of ${\mathcal E}$-integral curves 
$c : I \to {\mathcal F}_{\Lag}(\R^{2n+2})$ 
such that the osculating-framed contact-integral curve 
$\pi_1\circ c : I \to P(V^{2n+2})$ is of type $\AAA = (a_1, a_2, \dots, a_{2n+1})$ 
is not empty 
if and only if 
$$
a_{n+j} = a_{n+1} + a_n - a_{n+1-j}, \quad (2 \leq j \leq n+1), 
$$
and then its codimension in the jet space of ${\mathcal E}$-integral curves is given by 
$a_{n+1} - (n+1)$. 
\ent

\Proof
To show Theorem \ref{codim-E}, first 
we give systems of projective coordinates on ${\mathcal F}_{\Lag}(V)$. 
For the case $n = 1$, refer the paper \cite{IMT}. 

We fix a flag 
${\mathbf V}_0 = (V_{10}, V_{20}, \dots, V_{n+1\,0}) \in {\mathcal F}_{\Lag}(V)$. 
Then we take the open set $U \subset {\mathcal F}_{\Lag}(V)$ defined by 
$$
U := 
\left\{ 
(V_1, V_2, \dots, V_{n+1}) \in {\mathcal F}_{\Lag}(V) \ 
\right\vert  
\left.V_1\cap V_{10}^s = \{ 0\}, V_2\cap V_{20}^s = \{ 0\}, \dots, V_{n+1}\cap V_{n+1\,0}^s = \{ 0\} \right\}.
$$
Take ${\mathbf V}_1 = (V_{11}, V_{21}, \dots, V_{n+1\, 1}) \in U$. 
Then we have the decomposition
$
V  =  V_{n+1\, 1} \oplus V_{n+1\, 0} 
$
into Lagrangian subspaces, and the decomposition
$$
\begin{array}{rcl}
V_{n+1\, 1} & = & 
V_{11}\oplus(V_{21}\cap V_{10}^s)\oplus(V_{31}\cap V_{20}^s)\oplus\cdots\oplus(V_{n+1\, 1}\cap V_{n\, 0}^s), 
\vspace{0.2truecm}
\\
V_{n+1\, 0} & = & 
V_{10}\oplus(V_{20}\cap V_{11}^s)\oplus(V_{30}\cap V_{21}^s)\oplus\cdots\oplus(V_{n+1\, 0}\cap V_{n\, 1}^s), 
\end{array}
$$
of each Lagrangian subspace into one-dimensional subspaces. 
Take non-zero vectors $e_0 \in V_{11}$, $e_i \in V_{i+1\, 1}\cap V_{i\, 0}^s, (1 \leq i \leq n)$, $f_0 \in V_{10}$ and 
$f_i \in V_{i+1\, 0}\cap V_{i\, 1}^s, (1 \leq i \leq n)$, to get a symplectic basis 
$(e_0, e_1, \dots, e_n; f_0, f_1, \dots, f_n)$ of $V$. 

Then, for each ${\mathbf V} = (V_{1}, V_{2}, \dots, V_{n+1}) \in U_{{\mathbf V}_0}$, 
$V_{n+1}$ has a basis $(v_0, v_1, \dots, v_n)$ uniquely expressed as 
$$
v_i  =  e_i + {\displaystyle \sum_{j=0}^n} x_j^{\ i}f_j, \quad (0 \leq i \leq n), 
$$
for some $( x_j^{\ i})_{0 \leq i, j \leq n}$. 
Since $V_{n+1}$ is a Lagrangian subspace of $V$, 
we have that $x_j^{\ i} = x_i^{\ j}, 0 \leq i, j \leq n$. 
Then there exist uniquely $\lambda_i^{\ k}, (1 \leq k \leq i \leq n)$, 
such that 
$$
w_k = v_{k-1} + {\displaystyle \sum_{i=k}^n} \lambda_i^{\ k}v_i, \quad (1 \leq k \leq n+1), 
$$
form a basis of $V_{n+1}$ such 
that $V_k = \left\langle w_1, \dots, w_k\right\rangle_{\R}, (1 \leq k \leq n+1)$. 
Then actually we have 
$$
w_k = e_{k-1} + {\displaystyle \sum_{i=k}^n} \lambda_i^{\ k}e_i  
+ {\displaystyle \sum_{j=0}^n}\left( x_j^{\ k-1} + {\displaystyle \sum_{i=k}^n} \lambda_i^{\ k}x_j^{\ i}\right)f_j, 
\quad (1 \leq k \leq n+1). 
$$
Thus, given ${\mathbf V}_0, {\mathbf V}_1 \in {\mathcal F}_{\Lag}(V)$, 
we have a chart $U \to \R^{(n+1)^2}$ of ${\mathcal F}_{\Lag}(V)$, 
given by the symmetric matrix $( x_j^{\ i})_{0 \leq i, j \leq n}$ and $\lambda_i^{\ k}, (1 \leq k \leq i \leq n)$. 
From another choice of ${\mathbf V}_0, {\mathbf V}_1 \in {\mathcal F}_{\Lag}(V)$, we have 
another chart with fractional linear transition functions. 

The projection 
$\pi_1 : {\mathcal F}_{\Lag}(V) \to P(V)$ is expressed by 
$$
\left( x_j^{\ i}, \lambda_i^{\ k} \right) \mapsto 
\left[ 1: \lambda_1^{\ 1} : \cdots : \lambda_1^{\ n} : 
x_0^{\ 0} + {\displaystyle \sum_{i=1}^n} \lambda_i^{\ 1}x_0^{\ i} : \cdots : 
x_n^{\ 0} + {\displaystyle \sum_{i=1}^n} \lambda_i^{\ 1}x_n^{\ i}\right]. 
$$
We set 
$X_j^{\ k} := x_j^{\ k} + {\displaystyle \sum_{i=k+1}^n} \lambda_i^{\ k+1}x_j^{\ i}, 
(0 \leq j \leq n, 0 \leq k \leq n)$. 
Then the differential system ${\mathcal E}$ is locally given by 
$$
\left\{
\begin{array}{l}
d\lambda_i^{\ k} - \lambda_i^{\ k+1}d\lambda_k^{\ k} = 0, \quad 1 \leq k \leq n, k+1 \leq i \leq n, 
\vspace{0.2truecm}
\\
dX_j^{\ k-1} - X_j^{\ k} d\lambda_k^{\ k} = 0, \quad 1 \leq k \leq n, 0 \leq j \leq n. 
\end{array}
\right.
$$

We see that each ${\mathcal E}$-integral curve $\Gamma$ 
is obtained from the components $\lambda_k^{\ k}, 1 \leq k \leq n$, and 
the $x_n^{\ n}$-component, by iterative integrations. 

The type $(a_1, a_2, \dots, a_n, a_{n+1}, a_{n+2}, \dots, a_{2n+1})$ of 
$\gamma = \pi_1\circ\Gamma$ is expressed in terms of 
$$
u_k := \ord(\lambda_k^{\ k}), \ 1 \leq k \leq n, \quad v := \ord(x_n^{\ n})
$$
by 
$$
\begin{array}{rcl}
a_i & = & u_1 + u_2 + \cdots + u_i, \ (1 \leq i \leq n)
\\
a_{n+1} & = & u_1 + u_2 + \cdots + u_n + v, 
\\ 
a_{n+1+j} & = & u_1 + u_2 + \cdots + 2u_{n-j+1} + \cdots + 2u_n + v, \ (1 \leq j \leq n), 
\end{array}
$$
Let $\AAA = (a_1, \dots, a_n, a_{n+1}, a_{n+2}, \dots, a_{2n}, a_{2n+1})$ be a strictly 
increasing sequence of positive integers. Then 
The above system of equations has an integer solution $(u_1, \dots, u_n, v)$ if and only if 
$a_{n+1+i} - a_{n+i} = a_n - a_{n-i}$. 
If the non-empty condition is fulfilled, 
then the codimension of the set 
$$
\Sigma(\AAA) = \{ j^r\Gamma(t_0) \mid 
\Gamma : (I, t_0) \to {\mathcal F}_{\Lag}(V) {\mbox{\rm \ is }} {\mathcal E}{\mbox{\rm{-integral}}}, 
\ \type(\pi_1\circ\Gamma) = \AAA\} 
$$
in $J^r_{\mathcal E}(I, {\mathcal F}_{\Lag}(V))$ is calculated by 
$$
a_1 - 1 + (a_2 - a_1 - 1) + \cdots + (a_{n+1} - a_n - 1) = a_{n+1} - (n+1). 
$$
\QED

\

By Theorem \ref{codim-E} and by the transversality theorem for ${\mathcal E}$-integral curves, we have the following result: 
We separate cases into three groups from 
the classification viewpoint of singularities. 

\begin{figure}
\begin{center}
 \includegraphics[width=10truecm, height=2.5truecm, clip, bb=13.6421 681.004 574.086 830.14 ]{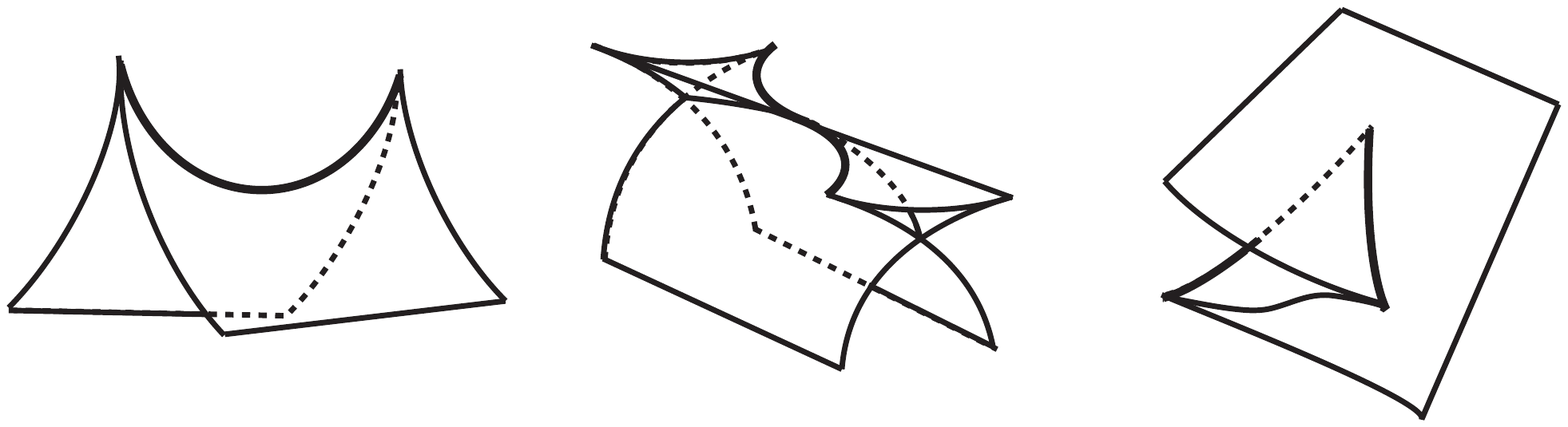}  
\caption{cuspidal edge, Mond surface and generic folded pleat in $\R^3$. }
\label{framed-contact}
   \end{center}
\end{figure}

\bet
\label{MIT}
{\rm (\cite{IMT})}
Let $2n+1 = 3$. For a generic ${\mathcal E}$-integral curve $c : I \to {\mathcal F}_{\Lag}(\R^4)$ 
in $C^\infty$-topology,  
the type $\AAA$ of $\pi_1\circ c$ at any point $t \in I$ is given by 
$$
\AAA = (1, 2, 3), (1, 3, 4), (2, 3, 5).
$$
The tangent varieties to the osculating-framed Legendre curve 
$\gamma = \pi_1\circ c : I \to P(V) \cong \R P^{3}$ is locally diffeomorphic to 
the cuspidal edge, to the Mond surface or to the generic folded pleat {\rm (Figure \ref{framed-contact})}. 
\ent

\ber
{\rm  
In the above Theorem \ref{MIT}, 
the type of the curve $\gamma$ is restricted to 
$(1, 2, 3), (1, 3, 4)$ or $(2, 3, 5)$. The local diffeomorphism class of the tangent variety $\Tan(\gamma)$ 
is determined if $\type(\gamma) = (1, 2, 3)$ or $(1, 3, 4)$, but it is not determined 
if $\type(\gamma) = (2, 3, 5)$ and there are exactly two diffeomorphism classes, 
{\it generic} one and {\it non-generic} one. 
}
\enr

Note that we have obtained in \cite{IMT} also the generic classification of singularities of tangent varieties 
to $\pi_2\circ c : I \to \LG(V)$ 
in Lagrangian Grassmannian. 

In the higher codimensional case, we have: 

\bet
\label{osculating framed contact-integral-geq7}
Let $2n + 1 \geq 7$. For a generic ${\mathcal E}$-integral curve $c : I \to {\mathcal F}_{\Lag}(\R^{2n+2})$ 
in $C^\infty$-topology,  the type of osculating-framed contact-integral curve 
$\gamma = \pi_1\circ c : I \to P(V) \cong \R P^{2n+1}$ at each point of $I$ is given by one of 
$$
\begin{array}{ccc}
\AAA & = & (1, 2, 3, 4, \dots, n, n+1, n+2, \dots, 2n+1), 
\\
& & (1, 2, 3, 4, \dots, n, \ \ n+2, n+3, \dots, 2n+2), 
\\
& & \cdots\cdots\cdots
\\
& & (1, 2, 4, 5, \dots, n+1, n+2, n+3, \dots, 2n+2), 
\\
& & (1, 3, 4, 5, \dots, n+1, n+2, n+3, \dots, 2n+2), 
\\
& & (2, 3, 4, 5, \dots, n+1, n+2, n+3, \dots, 2n+2). 
\end{array}
$$
Moreover the tangent variety $\Tan(\gamma)$ to the osculating-framed contact-integral curve 
$\gamma$ is locally diffeomorphic to 
the cuspidal edge, the open folded umbrella, the open Mond surface, or to the open swallowtail. 
\ent

We should be careful in the low codimensional case: 

\bet
\label{osculating framed contact-integral-=5}
Let $2n + 1 = 5$. For a generic ${\mathcal E}$-integral curve $c : I \to {\mathcal F}_{\Lag}(\R^6)$ 
in $C^\infty$-topology,  the type of osculating-framed contact-integral curve 
$\gamma = \pi_1\circ c : I \to P(V) \cong \R P^{5}$ at each point of $I$ is given by one of 
$$
(1, 2, 3, 4, 5), \ (1, 2, 4, 5, 6), \ (1, 3, 4, 6, 7), \ (2, 3, 4, 5, 7). 
$$
Moreover the tangent variety $\Tan(\gamma)$ to the osculating-framed contact-integral curve 
$\gamma$ is locally diffeomorphic to 
the cuspidal edge, the open folded umbrella, 
the unfurled Mond surface, or to the open swallowtail. 
\ent

\no
{\it Proofs of Theorems \ref{osculating framed contact-integral-geq7}, \ref{osculating framed contact-integral-=5}:}
By the transversality theorem, we reduce the list in each case from Theorem \ref{codim-E}. 
In each case, we have the uniqueness of the diffeomorphism class of tangent varieties by Theorem \ref{normal forms}, 
except for the case $\AAA = (1, 3, 4, 6, 7)$. 
For the case $\AAA = (1, 3, 4, 6, 7)$, we use
Theorem \ref{1,3,4,6}. 
\QED

\

It is natural to consider the generic classification of tangent varieties to contact-integral curves 
$I \to P(V) = \R P^{2n+1}$. 
Here, we give just the result on non-framed three dimensional case ($n = 1$): 

\bep
For a generic contact-integral curve 
$\gamma : I \to P(V^4) \cong \R P^{3}$, and for any $t_0 \in I$, 
the type of $\gamma$ at $t_0$ is equal to $(1, 2, 3)$ or to $(1, 3, 4)$ and 
the tangent variety 
$\Tan(\gamma)$ of $\gamma$ is locally diffeomorphic to 
the cuspidal edge or to the Mond surface. 
\enp

\Proof
Take the local coordinates $\lambda, \mu, \nu$ of $P(V)$ such that 
the contact structure is given by $d\mu = \nu d\lambda - \lambda d\nu$. 
We express $\gamma(t) = (\lambda(t), \mu(t), \nu(t))$. 
Since $\gamma$ is contact-integral, we have that 
$\mu'(t) = \nu(t)\lambda'(t) - \lambda(t)\nu'(t)$. 
Therefore $\mu''(t) = \nu(t)\lambda''(t) - \lambda(t)\nu''(t)$ and 
$$
\mu'''(t) = \nu'(t)\lambda''(t) + \nu(t)\lambda'''(t) - \lambda'(t)\nu''(t) - \lambda(t)\nu'''(t). 
$$
Then 
$$
\det\left(
\begin{array}{ccc}
\lambda' & \mu' & \nu' 
\\
\lambda'' & \mu'' & \nu'' 
\\
\lambda''' & \mu''' & \nu'''
\end{array}
\right) 
= (\lambda'\nu'' - \lambda''\nu')^2. 
$$
Therefore, if $\type(\lambda(t), \nu(t)) = (1, 2)$, then $\type(\gamma(t)) = (1, 2, 3)$. 
Moreover we have 
$$
\mu'''' = 2\nu'\lambda''' + \nu\lambda'''' - 2\lambda'\nu''' - \lambda\nu''''. 
$$
Then
$$
\rank
\left(
\begin{array}{ccc}
\lambda' & \mu' & \nu' 
\\
\lambda'' & \mu'' & \nu'' 
\\
\lambda''' & \mu''' & \nu'''
\\
\lambda'''' & \mu'''' & \nu''''
\end{array}
\right) 
= 
\rank 
\left(
\begin{array}{ccc}
\lambda' & \nu' & 0
\\
\lambda'' & \nu'' & 0
\\
\lambda''' & \nu''' & \lambda'\nu'' - \lambda''\nu'
\\
\lambda'''' & \nu'''' & \lambda'\nu''' - \lambda'''\nu'
\end{array}
\right).
$$
Therefore the rank of the above matrix is $3$ at $t$ if and only if 
$\lambda'\nu'' - \lambda''\nu' \not= 0$ or $\lambda'\nu''' - \lambda'''\nu' \not= 0$ at $t$. 
By the transversality theorem, we have that, for a generic $\gamma$ and for any $t_0 \in I$, 
(a) $\lambda'(t_0)\nu''(t_0) - \lambda''(t_0)\nu'(t_0) \not= 0$ or 
(b) $\lambda'(t_0)\nu''(t_0) - \lambda''(t_0)\nu'(t_0) = 0$ and 
$\lambda'(t_0)\nu'''(t_0) - \lambda'''(t_0)\nu'(t_0) \not= 0$. 
In case (a), $\type(\gamma) = (1, 2, 3)$ at $t_0$. In case (b), $\type(\gamma) = (1, 3, 4)$ at $t_0$. 
Then, by Theorem \ref{normal forms}(1), we have the required result. 
\QED

\section{Singularities of tangent varieties to surfaces.}
\label{Singularities of tangent varieties to surfaces.}

First we observe that the tangent varieties to a generic smooth surface are not frontal. 

\bee
\label{Tangents to Veronese}
{\rm 
Let $V = \left\{ 
A = \left.\left( \begin{array}{ccc}
a_{11} & a_{12} & a_{13} \\
a_{12} & a_{22} & a_{23} \\
a_{13} & a_{23} & a_{33} 
\end{array}
\right) \ 
\right\vert \ 
3\times 3, \ 
{\mbox{\rm symmetric}}\right\}$, 
\vspace{0.2truecm}
\\
the vector space of 
quadratic form of variables $x, y, z$. Then $\dim(V) = 6$. 
Let $S = P(\{ \rank(A) = 1 \}) \subset P(V) \cong \R P^5$ be the Veronese surface. Then we see that 
the tangent variety consists of the projection of the locus of semi-indefinite matrices of rank $2$ and $S$. 
Note that the secant variety ${\rm Sec}(S)$, 
the closure of the union of secants connecting any pair of points on $S$, 
consists of the projection of the locus of matrices of rank $\leq 2$ : 
$$
\begin{array}{c}
\Tan(S) =  S \cup P(\{ \rank(A) = 2, \ {\mbox{\rm semi-indefinite}}\})
\\
\subsetneq 
{\mbox{\rm Sec}}(S) = P(\{ \rank(A) \leq 2 \}) 
\subsetneq P(V). 
\end{array}
$$

\begin{figure}
\begin{center}
 \includegraphics[width=6truecm, height=2.5truecm, clip, bb=58.6665 698.487 271.868 791.056]{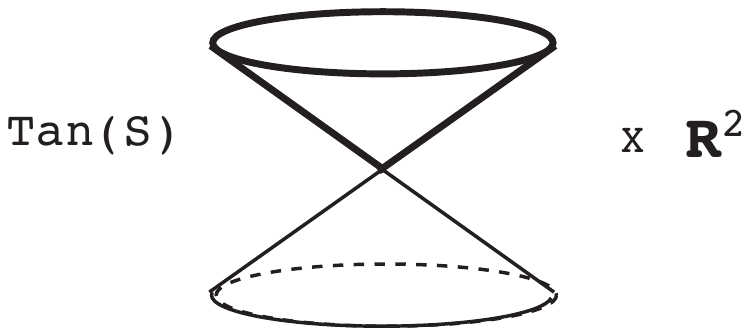}  
   \end{center}
   \caption{Tangent variety of Veronese surface. }
   \label{Veronese}
\end{figure}
}
\ene
See Figure \ref{Veronese}. 
The tangent variety $\Tan(S)$ is not frontal. 
Note that, even if $S$ is algebraic, $\Tan(S)$ is semi-algebraic in general over the real numbers.
For a generic surface $S \in \R P^5$, tangent varieties $\Tan(S)$ are perturbed into 
a non-frontal hypersurface. 

%



Therefore the tangent variety $\Tan(S)$ to a generic surface $S \subset \R P^5$ 
is never frontal. 

\

Let $V$ be a $(N+3)$-dimensional vector space. Let us consider a flag manifold 
$$
{\mathcal F} = {\mathcal F}_{1,3}(V) := 
\{ V_1 \subset V_3 \subset V \} \cong \Gr(2, T(P(V))), 
$$
${\mathcal F}_{1,3}(V) = 3N + 2$, with local coordinates $x_1, x_2, y_1, \dots, y_N, p_1, \dots, p_N, q_1, \dots, q_N$. 
The canonical differential system ${\mathcal T} = {\mathcal C} = {\mathcal C}_{1,3}$ is given by 
$
dy_i = p_i dx_1 + q_i dx_2, \quad (1 \leq i \leq N). 
$
A frontal map-germ $f : (\R^2, 0) \to P(V) = \R P^{N+2}$ lifts to a ${\mathcal C}_{1,3}$-integral map-germs, therefore 
$f$ is an openings of $g = (x_1\circ f, x_2\circ f) : (\R^2, 0) \to \R^2$ with the dense set of regular points. 

Thus it is possible to study the 
singularities of tangent varieties to frontal surfaces can be studied as the singularity theory on ${\mathcal C}_{1,3}$-integral mappings.
The general studies from this viewpoint are left to a forthcoming paper. 

Now, let us consider another type of flag manifold: 
$
{\mathcal F}_{1,3,5}(V) = \{ V_1 \subset V_3 \subset V_5 \subset V \}. 
$
and the canonical system ${\mathcal N} = {{\mathcal C}_{1,3,5}} \subset T({\mathcal F}_{1,3,5}(V))$ defined by 
$$
v \in {{\mathcal C}_{1,3,5}}_{(V_1, V_3, V_5)} 
\Longleftrightarrow 
\pi_{i*}(v) \in T(\Gr(i, V_{i+2})) (\subset T(\Gr(i, \R^6)), \ 
i = 1, 3. 
$$ 
If $N = 3$, 
then $\dim({\mathcal F}_{1,3,5}(\R^6)) = 13$ and $\rank({{\mathcal C}_{1,3,5}}) = 8$. 
In fact, ${\mathcal N}$ is given by 
$$
\left\{
\begin{array}{ccc}
dx_3^{\ 0} & = & x_3^{\ 1} dx_1^{\ 0} + x_3^{\ 2} dx_2^{\ 0}
\\
dx_4^{\ 0} & = & x_4^{\ 1} dx_1^{\ 0} + x_4^{\ 2} dx_2^{\ 0}
\\
dx_5^{\ 0} & = & x_5^{\ 1} dx_1^{\ 0} + x_5^{\ 2} dx_2^{\ 0}
\\
dx_5^{\ 1} & = & x_5^{\ 3} dx_3^{\ 1} + x_5^{\ 4} dx_4^{\ 1}
\\
dx_5^{\ 2} & = & x_5^{\ 3} dx_3^{\ 2} + x_5^{\ 4} dx_4^{\ 2} 
\end{array}
\right.
$$
for a system of projective local coordinates 
$$
x_1^{\ 0}, x_2^{\ 0}, x_3^{\ 0}, x_4^{\ 0}, x_5^{\ 0}, x_3^{\ 1}, x_4^{\ 1}, x_5^{\ 1}, 
x_3^{\ 2}, x_4^{\ 2}, x_5^{\ 2}, x_5^{\ 3}, x_5^{\ 4}
$$ 
of ${\mathcal F}_{1,3,5}(V^6)$. 

\bep
\label{1,3,5}
Let $f : (\R^2, 0) \to P(V^{N+3})$ be a frontal map-germ. 
Suppose that the regular locus of the tangent map $\Tan(f) : (\R^4, 0) \to P(V)$ is dense. 
Then $\Tan(f)$ is frontal if and only if 
$f$ is the projection of a ${{\mathcal C}_{1,3,5}}$-integral map by 
$\pi_1 : {\mathcal F}_{1,3,5}(V) \to P(V)$. 
\enp

\Proof
Suppose $\Tan(f)$ is frontal and $g : (\R^4, 0) \to 
\Gr(4, T(P(V))) = {\mathcal F}_{1,5}(V)$ is the Grassmannian lifting of $\Tan(f)$. 
Then $g\vert_{\R^2\times 0}$ lifts a ${{\mathcal C}_{1,3,5}}$-integral map 
$F : (\R^2, 0) \to {\mathcal F}_{1,3,5}(V)$ and $\pi_1\circ F = f$. 
Conversely if $\pi_1\circ F = f$ for a ${{\mathcal C}_{1,3,5}}$-integral map $F$, 
then $\Tan(f)$ lifts to $G : (\R^4, 0) \to {\mathcal F}_{1,3,5}(V)$ by 
$G(s_1, s_2, t_1, t_2) = F(0, 0, t_1, t_2)$. 
\QED

\

Let $V^6$ be a symplectic vector space. 
Let us consider the canonical contact structure on $P(V) = \R P^5$. 
Let $S \subset \R P^5$ be a Legendre surface. 
Then $S$ lifts to a ${{\mathcal C}_{1,3,5}}$-integral surface. Therefore, by Theorem \ref{1,3,5}, we have: 

\bec
Let $i : (\R^2, 0) \to \R P^5$ be a Legendre immersion-germ. Suppose the regular locus 
$\Reg(\Tan(i))$ of the tangent variety is dense in $(\R^2, 0)$. Then 
the tangent variety $\Tan(i) : (\R^2, 0) \to \R P^5$ is a frontal. 
\enc



\bef
{\rm 
A point $p$ of a Legendre surface $S$ in $\R P^5$ is called an {\it ordinary point} if 
there exists a local projective-contact coordinates $x_1, x_2, x_3, x_4, x_5$ and 
a $C^\infty$ local coordinates $(u, v)$ of $S$ centred $p$ such that 
locally $S$ is given by 
$$
\left\{
\begin{array}{lcl}
x_1 & = & u, 
\\
x_2 & = & v, 
\\
x_3 & = & \frac{1}{2}au^2 + buv + \frac{1}{2}cv^2 \ + \ {\mbox{\rm \footnotesize higher order terms}}, 
\vspace{0.2truecm}
\\
x_4 & = & \frac{1}{2}bu^2 + cuv  + \frac{1}{2}ev^2 \ + \ {\mbox{\rm \footnotesize higher order terms}}, 
\vspace{0.2truecm}
\\
x_5 & = & -(\frac{1}{6}au^3 + \frac{1}{2}bu^2v + \frac{1}{2}cuv^2 + \frac{1}{6}ev^3) \ + \ {\mbox{\rm \footnotesize higher order terms}}, 
\end{array}
\right.
$$
with 
$$
{\mathcal D} = \{ dx_5 - x_1dx_3 - x_2dx_4 + x_3dx_1 + x_4dx_2 = 0\}, 
$$
and 
$$
\rank \left(\begin{array}{ccc}
a & b & c \\
b & c & e 
\end{array}
\right) = 2. 
$$
An ordinary point $p$ is called 
{\it hyperbolic} (resp. {\it elliptic}, {\it parabolic}), if moreover 
$$
H := 4(ac - b^2)(be - c^2) - (ae - bc)^2
$$
is negative (resp. positive, zero). 
}
\enf

Note that the set of hyperbolic (resp. elliptic) ordinary points is an open subset in $S$. 
Then we have the following fundamental result: 

\bet
\label{tangent to ordinary}
The tangent variety $\Tan(S)$ to a Legendre surface $S$ in $\R P^5$ at 
a hyperbolic ordinary point {\rm(}resp. an elliptic ordinary point{\rm)} is 
locally diffeomorphic to 
$(D_4^+$-singularity in $\R^3)\times \R^2$ {\rm(}resp. 
$(D_4^-$-singularity in $\R^3)\times \R^2${\rm)}  in $\R^5$. 
\ent

\begin{figure}[htbp]
\vspace{0.5truecm}
\begin{center}
 \includegraphics[width=12truecm, height=2.5truecm, clip, bb=43.666 741.525 443.201 828.081]{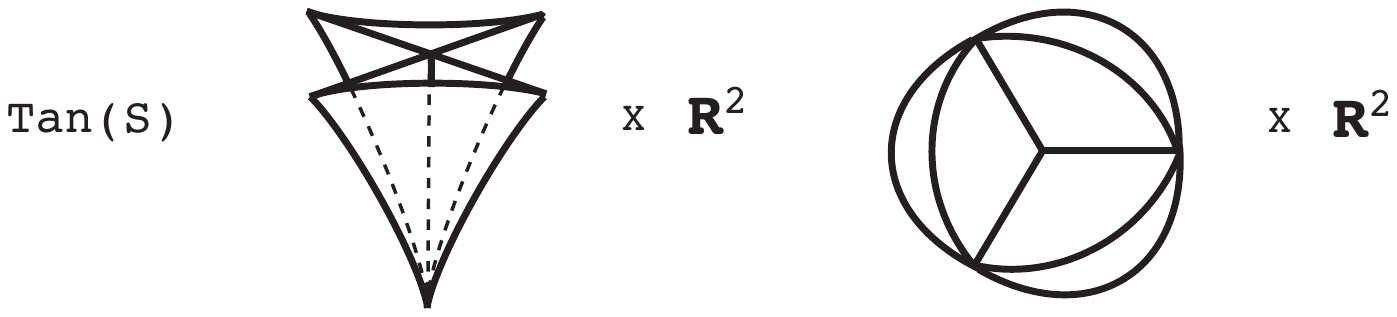}  
   \end{center}
   \vspace{-0.5truecm}
   \caption{Tangent varieties along hyperbolic and elliptic ordinary points on a surface in $\R P^5$.}
   \label{D_4}
\end{figure}

In \cite{Saji}, a simple criterion on $D_4$ has been found by Saji. 
The $D_4^{\pm}$-singularity in $\R^3$ is given by the map-germ $(\R^2, 0) \to (\R^3, 0)$
$$
(u, v) \mapsto (uv, u^2 \pm 3v^2, u^2v \pm v^3). 
$$

\bet
\label{Saji's criterion}
{\rm (\cite{Saji})}
Let $f : (\R^2, 0) \to (\R^3, 0)$ be a front and $(f, \nu) : (\R^2, 0) 
\to \R^3\times S^2$ a Legendre lift of $f$. Then 
$f$ is diffeomorphic to $D_4^+$ {\rm(}resp. $D_4^-${\rm)} 
if and only if $f$ is of rank zero at $0$ and 
the Hessian determinant of 
$$
\lambda(u, v) := \det\left(\dfrac{\pa f}{\pa u}(u, v), \dfrac{\pa f}{\pa v}(u, v), \nu(u, v) \right)
$$
at $(0, 0)$ is negative {\rm(}resp. positive{\rm)}. 
\ent

Note that $D_4$-singularity is not a generic singularity of wave-fronts in $\R^3$, 
but is a generic singularity of wave-fronts in $\R^4$. 
The criterion for $D_4$-singularities in $\R^4$ is also given in \cite{Saji}. 
Moreover we remark that Saji's criterion is valid also for the case with parameters and it 
characterises the trivial deformation of $D_4$-singularity. In fact the same line of proof in \cite{Saji} 
works as well for the case with parameters: 

\bet
\label{Saji's criterion2}
Let $F = (f_t)_{t \in (\R^r, 0)} : (\R^2\times\R^r, 0) \to (\R^3, 0)$ 
be a family of fronts and $(F, N) = (f_t, \nu_t) : (\R^2\times\R^r, 0) 
\to \R^3\times S^2$ a family of Legendre lifts of $F$. Then 
$F$ is diffeomorphic to the trivial deformation of $D_4^+$ {\rm(}resp. $D_4^-${\rm)} 
if and only if $f_t$ is of rank zero at $0$ and 
the Hessian determinant of 
$$
\lambda(u, v, t) := \det\left(\dfrac{\pa f_t}{\pa u}(u, v), \dfrac{\pa f_t}{\pa v}(u, v), \nu_t(u, v) \right)
$$
with respect to $(u, v)$ at $(0, 0, t)$ is negative {\rm(}resp. positive{\rm)}, for any $t \in (\R^r, 0)$. 
\ent

\

\noindent
{\it Proof of Theorem \ref{tangent to ordinary}:} 
Let $x_1 = u, x_2 = v$, 
$$
\begin{array}{lcl}
x_3 & = & \frac{1}{2}au^2 + buv + \frac{1}{2}cv^2 + \varphi(u, v), 
\vspace{0.2truecm}
\\
x_4 & = & \frac{1}{2}bu^2 + cuv  + \frac{1}{2}ev^2 + \psi(u, v), 
\vspace{0.2truecm}
\\
x_5 & = & -(\frac{1}{6}au^3 + \frac{1}{2}bu^2v + \frac{1}{2}cuv^2 + \frac{1}{6}ev^3) + \rho(u, v),  
\end{array}
$$
$\ord(\varphi) \geq 3, \ord(\psi) \geq 3$, and 
$$
\rho_u = u\varphi_u + v\psi_u - \varphi, \quad 
\rho_u = u\varphi_u + v\psi_u - \psi. 
$$
As an integrability condition, we have that $\varphi_v = \psi_u$. 
The tangent map of $S$ is given by $x_1 = u + s,  x_2  = v + t$, 
$$
\begin{array}{lcl}
x_3 & = & \frac{1}{2}au^2 + buv + \frac{1}{2}cv^2 + \varphi + 
s\left(au + bv + \varphi_u\right) 
+ t\left(bu + cv + \varphi_v \right), 
\vspace{0.1truecm}
\\
x_4 & = & \frac{1}{2}bu^2 + cuv  + \frac{1}{2}ev^2 + \psi + s\left(bu + cv + \psi_u\right) 
+ t\left(cu + ev + \psi_v \right), 
\vspace{0.1truecm}
\\
x_5 & = & -(\frac{1}{6}au^3 + \frac{1}{2}bu^2v + \frac{1}{2}cuv^2 + \frac{1}{6}ev^3) + \rho 
\\
& & \quad 
+ s\left(-\frac{1}{2}au^2 - buv - \frac{1}{2}cv^2 + \rho_u\right) 
+ t\left(-\frac{1}{2}bu^2 - cuv - \frac{1}{2}ev^2 + \rho_v\right). 
\end{array}
$$
Take the transversal slice $s = -u, t = -v$. Then we have map-germ $g : (\R^2, 0) \to (\R^3, 0)$, 
$$
\begin{array}{lcl}
g_1(u, v) & = & -\frac{1}{2}au^2 - buv - \frac{1}{2}cv^2 + \varphi - u\varphi_u - v\varphi_v, 
\vspace{0.1truecm}
\\
g_2(u, v) & = & -\frac{1}{2}bu^2 - cuv - \frac{1}{2}ev^2 + \psi - u\psi_u - v\psi_v, 
\vspace{0.1truecm}
\\
g_3(u, v) & = & \frac{1}{3}au^3 + bu^2v + cuv^2 + \frac{1}{3}ev^3 + \rho - u\rho_u - v\rho_v. 
\end{array}
$$
We show that $g$ is diffeomorphic to $D_4$-singularity, by using Saji's criterion (Theorem \ref{Saji's criterion}). 

First, we have $dg_3 = - u dg_1 - v dg_2$. Therefore $g$ is a front and we can take 
$\nu = \frac{1}{\sqrt{u^2+v^2+1}}\left(u, v, 1\right)$. 
Second, we see $f$ is of rank zero. Third, 
$$
\begin{array}{rcl}
\lambda(u, v) & = & \det(g_u, g_v, \nu) = \det\left(
\begin{array}{ccc}
g_{1u} & g_{1v} & u 
\\
g_{2u} & g_{2v} & v 
\\
0 & 0 & \sqrt{u^2+v^2+1}
\end{array}
\right) 
\vspace{0.2truecm}
\\
& = & \sqrt{u^2+v^2+1}\left(g_{1u}g_{2v} - g_{1v}g_{2u}\right)
\end{array}
$$
The $2$-jet of $h := g_{1u}g_{2v} - g_{1v}g_{2u}$ at $0$ is given by
$$
j^2h(0) = (ac - b^2)u^2 + (ae - bc)uv + (be - c^2)v^2 \quad (\mod. \ {\mathfrak m}_2^3). 
$$
Therefore we have that 
the Hessian determinant of $\lambda$ at $0$ is given by 
$$
H = \det
\left(
\begin{array}{cc}
2(ac - b^2) & ae - bc 
\\
ae - bc & 2(be - c^2) 
\end{array}
\right)
$$ 
By Theorem \ref{Saji's criterion}, 
we see that $g$ is diffeomorphic to $D_4^{\pm}$ if and only if $\mp H > 0$. 
Moreover, we can show similarly that, regarding $S$ as the parameter space, 
the tangent map-germ is diffeomorphic to the trivial unfolding of $D_4$-singularity with two parameters, 
by using Theorem \ref{Saji's criterion2}. 
Hence we have Theorem \ref{tangent to ordinary}. 
\QED

\section{Tangent maps to frontal maps and open problems.}
\label{Tangent maps to frontal maps.}

Let $V$ be a $(N+2n)$-dimensional vector space with positive natural numbers $N, n$.  
Consider the flag manifolds:
$$
{\mathcal F}_{1,n+1,2n+1} = {\mathcal F}_{1,n+1,2n+1}(V) := 
\{ V_1 \subset V_{n+1} \subset V_{2n+1} \subset V\}, 
$$
with the canonical differential system ${\mathcal C}_{1,n+1,2n+1}$, and 
$$
{\mathcal F}_{1,n+1} = {\mathcal F}_{1,n+1}(V) := 
\{ V_1 \subset V_{n+1} \subset V\}, 
$$
with the canonical differential system ${\mathcal C}_{1,n+1}$. 
Note that ${\mathcal F}_{1,n+1}$ is identified with the Grassmannian bundle 
$\Gr(n, T(P(V)))$. 
Consider the canonical projections 
$$
{\mathcal F}_{1,n+1,2n+1} \stackrel{\Pi}{\longrightarrow} {\mathcal F}_{1,n+1} 
\stackrel{\pi}{\longrightarrow} {\mathcal F}_{1} = P(V) = \R P^{N+2n-1}. 
$$

Similarly to the proof of Proposition \ref{1,3,5}, we have 

\bep
\label{1,n+1,2n+1}
Let $f : (\R^n, 0) \to \R P^{N+2n-1}$ be a frontal map-germ. 
Suppose the regular locus $\Reg(\Tan(f))$ of 
the tangent map $\Tan(f) : (\R^{2n}, 0) \to \R P^{N+2n-1}$ is dense in $(\R^{2n}, 0)$. 
Then $\Tan(f)$ is frontal if and only if 
the Grassmannian lift $\widetilde{f} :  (\R^n, 0) \to {\mathcal F}_{1,n+1}$ of $f$ for $\pi$, 
lifts to a ${\mathcal C}_{1,n+1,2n+1}$-integral lift 
${\mathbf f} :  (\R^n, 0) \to {\mathcal F}_{1,n+1,2n+1}$ for $\Pi$. 
\enp

It is natural to proceed to consider the tangent varieties to Legendre submanifolds. 

Let $V$ be a $(2n+2)$-dimensional symplectic vector space. 
Consider the Lagrange (isotropic) flag manifold:
$$
{\mathcal F}_{\Lag} = {\mathcal F}_{\Lag}(V) := 
\{ V_1 \subset V_{n+1} \subset V \ \mid \ V_{n+1} {\mbox{\rm\ is Lagrange.}} \}, 
$$
with the canonical differential system ${\mathcal E} \subset T{\mathcal F}_{\Lag}$. 
In general we have 

\bec
\label{E-integral-general}
Let $g : (\R^n, 0) \to {\mathcal F}_{\Lag}$ be ${\mathcal E}$-integral 
and 
$\Tan(\pi_1\circ g) : (\R^{2n}, 0) \to P(V)$ the tangent map-germ of 
$\pi_1\circ g : (\R^n, 0) \to P(V)$. 
Suppose that $\Reg(\Tan(\pi_1\circ g))$ is dense in $(\R^n, 0)$.   
Then $\Tan(\pi_1\circ g)$ is frontal. 
\enc

\Proof
Note that ${\mathcal F}_{\Lag}$ is embedded in ${\mathcal F}_{1,n+1,2n+1}$ 
by $(V_1, V_{n+1}) \mapsto (V_1, V_{n+1}, V_1^s)$, where $V_1^s$ is the symplectic 
skew-orthogonal to $V_1$, and ${\mathcal E}$ is the restriction of 
${\mathcal C}_{1,n+1,2n+1}$. Therefore Proposition \ref{E-integral-general} 
follows from Proposition \ref{1,n+1,2n+1}. 

Here we give alternative direct proof. 
Since $f$ is Legendre,  $f = (\lambda, \mu, \nu)$ satisfies 
$d\mu = \sum_{i=1}^n \left(\nu_id\lambda_i - \lambda_id\nu_i\right)$. 
The tangent map-germ $\Tan(f) = (\Lambda, M, N)$ is given by 
$$
\left(
\begin{array}{c}
\Lambda \\
M \\
N 
\end{array}
\right)
= 
\left(
\begin{array}{c}
\lambda \\
\mu \\
\nu 
\end{array}
\right)
+ 
\sum_{j=1}^n 
s_j
\left(
\begin{array}{c}
\pa \lambda/\pa u_j \\
\pa \mu/\pa u_j \\
\pa \nu/\pa u_j
\end{array}
\right).
$$
Then we have 
$$
\begin{array}{rcl}
dM & = & d\mu + 
{\displaystyle 
\sum_{j=1}^n
}
\
 s_j d\left(\pa\mu/\pa u_j\right) + \sum_{j=1}^n \left(\pa\mu/\pa u_j\right) ds_j
\\ 
& = & 
d\mu + 
{\displaystyle \sum_{i=1}^n\sum_{j=1}^n
}
\ 
s_j\left(\nu_id\left(\pa \lambda_i/\pa u_j\right) 
- 
\lambda_id\left(\pa \nu_i/\pa u_j\right)\right)
\\
&  & \quad\ 
+ 
{\displaystyle 
\sum_{i=1}^n\sum_{j=1}^n
}
\left(\nu_i\left(\pa \lambda_i/\pa u_j\right) - 
\lambda_i\left(\pa \nu_i/\pa u_j\right)\right)ds_j
\\
& = & 
{\displaystyle \sum_{i=1}^n
}
 \left( \nu_id\Lambda_i - \lambda_idN_i\right). 
\end{array}
$$
Thus $M \in {\mathcal R}_{(\Lambda, N)}$ and 
$\Tan(f)$ is frontal. 
\QED

Then Corollary \ref{E-integral-general} implies 

\bec
Let $f : (\R^n, 0) \to P(V) = \R P^{2n+1}$ be a germ of Legendre immersion and 
$\Tan(f) : (\R^{2n}, 0) \to P(V)$ the tangent map-germ of $f$. 
Suppose that $\Reg(\Tan(f))$ is dense in $(\R^n, 0)$.   
Then $\Tan(f)$ is frontal. 
\enc

\


We conclude the paper by posing open generic classification problems, which remain 
to be solved first: 
\\
{\it Problem 1:} Classify the singularities of tangent varieties to generic contact-integral curves in 
$P(V^{2n+2}) \cong \R P^{2n+1}$ for a symplectic vector space $V$ of dimension $2n+2$, under diffeomorphisms and 
contactomorphisms. 
\\
{\it Problem 2:} Classify the singularities of tangent varieties to generic surfaces in $\R P^5$.  
It would be natural to relate singularities of tangent variety to the method of height function or hight family 
(cf. \cite{Wall1}\cite{MFR}). 
\\
{\it Problem 3:} Classify the singularities of tangent varieties to generic frontal surfaces (projections of 
generic ${{\mathcal C}_{1,3}}$-integral surfaces in ${\mathcal F}_{1,3}(\R^6)$) in $\R P^5$.  
\\
{\it Problem 4:} Classify the singularities of tangent varieties to projections in $\R P^5$ of generic 
${{\mathcal C}_{1,3,5}}$-integral surfaces in ${\mathcal F}_{1,3,5}(\R^6)$. 
\\
{\it Problem 5:} 
Classify the singularities of tangent varieties to Legendre surfaces in $\R P^5$ 
along parabolic ordinary points. 
Moreover classify the singularities of tangent varieties of generic Legendre surfaces
in $\R P^5$. (See \S \ref{Singularities of tangent varieties to surfaces.}.)

\end{document}